\documentclass{agtart_a}
\pdfoutput=1

\usepackage{pinlabel}

\title[The algebraic crossing number and the braid index of knots and links]{The algebraic crossing number and\\the braid index of knots and links}

\author{Keiko Kawamuro}
\givenname{Keiko}
\surname{Kawamuro}
\address{Math Department\\Rice University\\
6100 S Main St\\\newline
Houston TX 77005-1892\\USA}
\email{keiko.kawamuro@rice.edu}
\urladdr{}

\volumenumber{6}
\issuenumber{}
\publicationyear{2006}
\papernumber{81}
\startpage{2313}
\endpage{2350}

\doi{}
\MR{}
\Zbl{}


\received{28 April 2006}
\revised{}
\accepted{21 July 2006}
\published{8 December 2006}
\publishedonline{8 December 2006}
\proposed{}
\seconded{}
\corresponding{}
\editor{CPR}
\version{}

\arxivreference{} %%% not in arxiv

\let\xysavmatrix\xymatrix
\def\xymatrix{\disablesubscriptcorrection\xysavmatrix}
\AtBeginDocument{\let\bar\wbar\let\tilde\wtilde\let\hat\what}
\def\strutt{\vrule width0pt height 12pt depth 5pt}
\def\struttt{\vrule width0pt height 12pt depth 8pt}
\allowdisplaybreaks

\makeatletter
\def\cnewtheorem#1[#2]#3{\newtheorem{#1}{#3}[section]
\expandafter\let\csname c@#1\endcsname\c@theorem}
\makeatother

\newtheorem{theorem}{Theorem}[section]
\cnewtheorem{lemma}[theorem]{Lemma}
\cnewtheorem{sublemma}[theorem]{Sublemma}
\cnewtheorem{proposition}[theorem]{Proposition}
\cnewtheorem{corollary}[theorem]{Corollary}
\theoremstyle{remark}
\cnewtheorem{definition}[theorem]{Definition}
\cnewtheorem{conjecture}[theorem]{Conjecture}
\cnewtheorem{remark}[theorem]{Remark}
\cnewtheorem{example}[theorem]{Example}

\newcommand{\ol}{\overline}
\newcommand{\lra}{\longrightarrow}
\newcommand{\Lra}{\Longrightarrow}

\begin{document}

\begin{htmlabstract}
<p class="noindent"> It has been conjectured that the algebraic
crossing number of a link is uniquely determined in minimal braid
representation. This conjecture is true for many classes of knots and
links.  </p> <p class="noindent"> The Morton-Franks-Williams
inequality gives a lower bound for braid index. And sharpness of the
inequality on a knot type implies the truth of the conjecture for the
knot type.  </p> <p class="noindent"> We prove that there are
infinitely many examples of knots and links for which the inequality
is not sharp but the conjecture is still true.  We also show that if
the conjecture is true for K and L, then it is also true for the
(p,q)-cable of K and for the connect sum of K and L.  </p>
\end{htmlabstract}

\begin{abstract}

It has been conjectured that the algebraic crossing number of a link
is uniquely determined in minimal braid representation. This
conjecture is true for many classes of knots and links.

The Morton--Franks--Williams inequality gives a lower bound for
braid index. And sharpness of the inequality on a knot type
implies the truth of the conjecture for the knot type.

We prove that there are infinitely many examples of knots and links
for which the inequality is not sharp but the conjecture is still
true.  We also show that if the conjecture is true for ${\cal K}$ and
${\cal L},$ then it is also true for the $(p,q)$--cable of ${\cal K}$
and for the connect sum of ${\cal K}$ and ${\cal L}.$
\end{abstract}

\begin{asciiabstract}
It has been conjectured that the algebraic crossing number of a link
is uniquely determined in minimal braid representation. This
conjecture is true for many classes of knots and links.

The Morton-Franks-Williams inequality gives a lower bound for
braid index. And sharpness of the inequality on a knot type
implies the truth of the conjecture for the knot type.

We prove that there are infinitely many examples of knots and links
for which the inequality is not sharp but the conjecture is still
true.  We also show that if the conjecture is true for K and L, then
it is also true for the (p,q)-cable of K and for the connect sum of K
and L.
\end{asciiabstract}

\maketitle

\section{Introduction}

The {\em braid index} is one of the classical invariants of knots and links.
Any knot and link type is presented as a braid closure.
The braid index of a link type is the least number of braid strands
needed for that.

\begin{figure}[ht!]\small
\begin{center}
\labellist
\pinlabel {$K_+$} at 51 15
\pinlabel {$K_-$} at 197 15
\pinlabel {$K_0$} at 340 15
\endlabellist
\includegraphics[height=20mm]{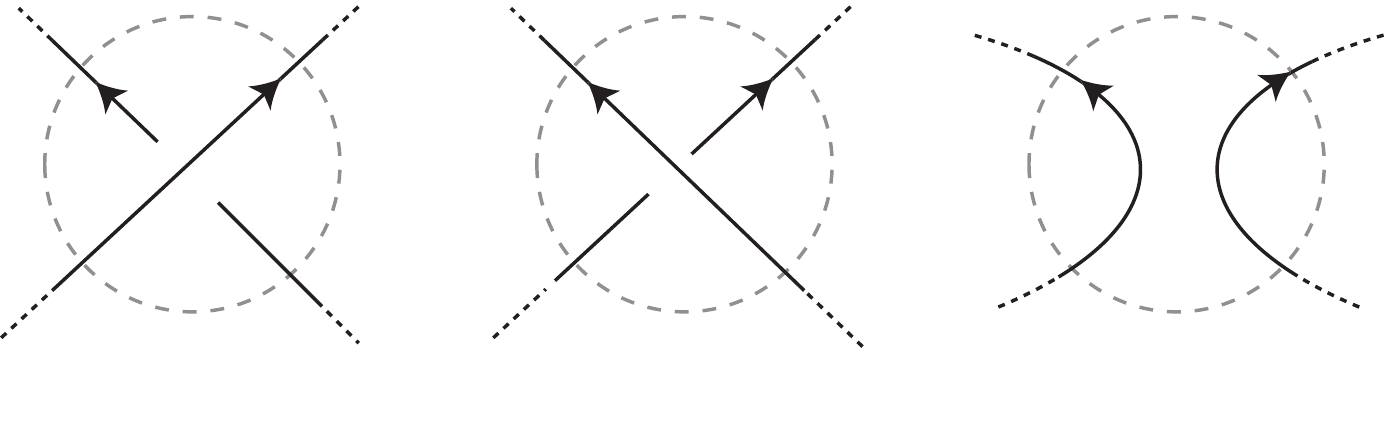}
\caption{Local views of $K_+, K_-, K_0$}\label{+-0}
\end{center}
\end{figure}

The {\em algebraic crossing number} (or writhe) is an integer 
associated to an oriented link diagram counting the crossings with
weight $+1$ (resp.\ $-1$) for a positive (resp.\ negative) crossing
as shown in the left (resp.\ middle) sketch of \fullref{+-0}.
Since it is changed under Reidemeister move I, it is {\em not} an
invariant of link types. However, it has been asked (see Jones
\cite[page 357]{Jones-1}  for example):

\medskip
\noindent\textbf{Question}\qua
{\sl Is the algebraic crossing number in a minimal braid
representation a link invariant?}

\medskip Here ``minimal'' means that the number of braid strands
of a link diagram is equal to the braid index of the link type.

It is known that the following links have unique algebraic crossing
numbers at minimal braid index: torus links, closed positive braids
with a full twist, including the Lorenz links (Franks and Williams
\cite{FW}), $2$--bridge links and alternating fibered links (Murasugi
\cite{Murasugi}) and links with braid index $\leq 3$ (Birman and
Menasco \cite{BM3}).

%We summarise \fullref{chap3} before \fullref{chap2}.

In \fullref{chap3} of this paper we approach the above question in
three ways.  The first way (\fullref{deficit-cable} and its
corollaries) is by studying the deficit of the
Morton--Franks--Williams (MFW) inequality (Morton \cite{Morton},
Franks and Williams \cite{FW}).  It is easy to see that sharpness of
the MFW--inequality implies the uniqueness of the algebraic crossing
number at minimal index.  Then how do we answer the question for links
on which the inequality is not sharp? In fact we provide infinitely
many examples of non-sharp links having unique algebraic crossing
numbers at minimal braid index.

The second way is by studying the behavior of the braid index and the
algebraic crossing number under the cabling operation. In \fullref{cable-thm} and \fullref{cable-thm-link},
we will prove that the uniqueness property is
preserved under cabling. Then we have \fullref{iterated-torus}
saying ``yes'' to the question for iterated torus knots.

The third way is by studying the connect sum operation. In \fullref{sum}, we will show that the uniqueness property is preserved
under taking the connect sum.

In \fullref{chap2}, we focus on non-sharpness of the MFW--inequality.

To state the MFW--inequality, let ${\cal K}$ be an oriented knot
type and let $K$ be a diagram of ${\cal K}$ on a plane. Focus
on one crossing of $K$ with sign $\varepsilon$. Denote
$K_\varepsilon := K$ and let $K_{-\varepsilon}$ (resp.\ $K_0$) be
the closed braid obtained from $K_\varepsilon$ by changing the the
crossing to the opposite sign $-\varepsilon$ (resp.\ resolving the
crossing),
see \fullref{+-0}. %

The \textit{HOMFLYPT polynomial} $P_{\cal K}(v,z)=P_{K}(v,z)$ satisfies the
following relations (for any choice of a crossing):
\begin{eqnarray}
\frac{1}{v} P_{K_{+}} - v P_{K_{-}} &=& z P_{K_{0}}.  \label{skein1} \\
P_{\mathrm{unknot}} &=&1.  \notag
\end{eqnarray}

Now we are ready to state the MFW--inequality.
\begin{theorem}[The Morton--Franks--Williams inequality \cite{Morton,FW}]
\label{MFW}
Let $d_{+}$ and $d_{-}$ be the maximal and minimal degrees of the
variable $v$ of $P_{\cal K}(v,z)$. If a knot type ${\cal K}$ has a
closed braid representative $K$ with braid index $b_K$ and
algebraic crossing number $c_K$, then we have
\begin{equation}
c_K - b_K +1 \leq d_{-}\leq d_{+}\leq c_K + b_K -1. \label{lower bound-1}%
\end{equation}%
As a corollary,
\begin{equation}
\frac{1}{2}(d_{+}-d_{-})+1 \leq b_K,  \label{lower bound}
\end{equation}%
giving a lower bound for the {\em braid index} $b_{\cal K}$ of
${\cal K}$.
\end{theorem}

In general, it is hard to determine the braid index. This
inequality was the first known result of a general nature relating
to the computation of braid index, and it appeared to be quite
effective.  Jones notes, in \cite{Jones-1}, that on all but five
knots, $9_{42}, 9_{49}, 10_{132}, 10_{150}, 10_{156}$ in the
standard knot table, up to crossing number $10$, the MFW
inequality is sharp. Furthermore it has been known that the
inequality is sharp on all torus links, closed positive $n$--braids
with a full twist \cite{FW}, $2$--bridge links and fibered
alternating links \cite{Murasugi}.

However, the MFW--inequality is not as strong as it appears to be
as above. In \fullref{deficit-thm} we give an infinite class
of prime links in which the
deficit $D_{\cal K}:= b_{\cal K}-\frac{1}{2}(d_+-d_-)-1$
of the MFW--inequality \eqref{lower bound}
can be arbitrarily large. And in \fullref{BM-thm} we see another infinite class of knots, including
$9_{42}, 9_{49}, 10_{132}, 10_{150}, 10_{156},$ on which the
inequality is not sharp.

Then we may ask ``why does non-sharpness occur?'' \fullref{thmA}
gives a sufficient condition for non-sharpness of the MFW
inequality. In fact all the examples in Theorems \ref{deficit-thm}
and \ref{BM-thm} satisfy this sufficient condition.

The idea of \fullref{thmA} is to find knots $K_\alpha$ of
known braid index $=b$ which have a distinguished crossing such
that, after changing that crossing to each of the other two
possibilities in \fullref{+-0}, giving knots or links $K_\beta$
and $K_\gamma$, it is revealed that $K_\beta$ and $K_\gamma$ each
has braid index $< b.$

Thanks to \fullref{thmA} one can visually observe the
``accumulation'' of deficits (for example under the connect sum
operation and other linking operation)
by looking only at the distinguished crossings which
contribute to deficits. See the proof of \fullref{deficit-thm}
for details.

\medskip
\textbf{Acknowledgment}\qua This paper is part of the author's PhD
thesis. She is grateful to her advisor, Joan Birman, for her
thoughtful advice and encouragement.  She also wishes to thank William
Menasco, who told her about the Birman--Menasco diagram and the
associated conjecture, when she visited SUNY Buffalo.  She appreciates
many helpful comments by Walter Neumann, Dylan Thurston, Ilya Kofman
and the referee and thanks Alexander Stoimenow for sending a preprint
She acknowledges partial support from NSF grants DMS-0405586 and
DMS-0306062.  Finally, she especially thanks Mikami Hirasawa, who
shared many creative ideas and results about fibered knots including
the definition and properties of the enhanced Milnor number.

%%%%%%%%%%%%%%%%%%%%%%%%%%%%%%%%%%%%%%%%%%%%%%%%%%%%%%%%%%%%%%%%%%%%%%%%%
\section{Non-sharpness of the Morton--Franks--Williams inequality}\label{chap2}

\subsection{Sufficient conditions for non-sharpness}
We define the deficit of MFW--inequality
(\fullref{deficit-def}) then give sufficient conditions
(\fullref{thmA}) for a closed braid on which the inequality is
not sharp.

Let $b_{\cal K}$ be the braid index of knot type ${\cal K}$,
that is the smallest integer $b_{\cal K}$ such that ${\cal K}$ can
be represented by a closed $b_{\cal K}$--braid. Let $b_K, c_K$ denote the
braid index and the algebraic crossing number of a braid
representative $K$ of ${\cal K}$.

\begin{definition}\label{deficit-def} Let
\begin{equation*}
D_{\cal K} : = b_{\cal K} - \frac{1}{2}(d_{+}-d_{-})-1
\end{equation*}%
be the difference of the numbers in $\eqref{lower bound},$ ie,
of the actual braid index and the lower bound for braid index.
Call $D_{\cal K}$ the {\em deficit} of the MFW--inequality for ${\cal K}$.
\end{definition}
If $D_{\cal K} =0,$ the MFW--inequality is sharp on ${\cal K}$. If
$K$ is a braid representative of ${\cal K}$ let $D_K^+ \ := \ (c_K
+ b_{K} -1)-d_{+}$ and $D_K^- \ := \ d_{-}-(c_K - b_{K} +1).$ When
$b_K = b_{\cal K},$ we have
\begin{equation}\label{D+D}
 D_{\cal K} =\frac{1}{2}(D_K^{+} + D_K^{-}).
\end{equation}%
Note that $D_K^{\pm}$ depends on the choice of braid
representative $K$, but the deficit $D_{\cal K}$ is independent
from the choice.

\begin{theorem}\label{thmA}
Assume that $K$ is a closed braid representative of ${\cal K}$ with
$b_K = b_{\cal K}$. Focus on one site of $K$ and construct
$K_+, K_-, K_0$ (one of the three must be $K$). Let $\alpha,
\beta, \gamma \in \{ +, -, 0 \}$ and assume that $\alpha, \beta,
\gamma$ are mutually distinct.
If $K_\alpha = K$ and positive
destabilization is applicable $p$--times
to each of $ K_\beta$ and $K_\gamma$, then
\begin{equation}
D_K^+ \geq  2p;    \label{less}
\end{equation}
and if $K_\alpha = K$ and negative destabilization is applicable $n$--times
to each of $ K_\beta$ and $K_\gamma$, then
\begin{equation}
D_K^- \geq  2n. \label{more}
\end{equation}
Therefore, by {\em \eqref{D+D}}, the MFW--inequality is not sharp on ${\cal K}$
if $p+n>0$.
\end{theorem}

Here is a lemma to prove \fullref{thmA}.
\begin{lemma}\label{lemma-for-thmA}
Let $K$ be a closed braid. Choose one crossing, and construct
$K_+, K_-, K_0$ (one of the three must be $K$). We have
\begin{eqnarray}
d_{+}(P_{K_{+}}) &\leq & \max \{d_{+}(P_{K_{-}})+2, \quad d_{+}(P_{K_{0}})+1\} \label{dplus+} \\
d_{+}(P_{K_{-}}) &\leq & \max \{d_{+}(P_{K_{+}})-2, \quad d_{+}(P_{K_{0}})-1\} \label{dplus-} \\
d_{+}(P_{K_{0}}) &\leq & \max \{d_{+}(P_{K_{+}})-1, \quad d_{+}(P_{K_{-}})+1\} \label{dplus0}
\end{eqnarray}%
and%
\begin{eqnarray*}
d_{-}(P_{K_{+}}) &\geq & \min \{d_{-}(P_{K_{-}})+2, \quad d_{-}(P_{K_{0}})+1\} \\
d_{-}(P_{K_{-}}) &\geq & \min \{d_{-}(P_{K_{+}})-2, \quad d_{-}(P_{K_{0}})-1\} \\
d_{-}(P_{K_{0}}) &\geq & \min \{d_{-}(P_{K_{+}})-1, \quad d_{-}(P_{K_{-}})+1\}. \end{eqnarray*}
\end{lemma}
\begin{proof}[Proof of \fullref{lemma-for-thmA}] 
By \eqref{skein1}, we have $P_{K_{+}}=v^{2}P_{K_{-}} + vz P_{K_{0}}.$
Thus, $d_+(P_{K_+})=d_+(v^2 P_{K_-} + vz P_{K_0}) \leq \max \{ d_+
( v^2 P_{K_-}), \ d_+ ( vz P_{K_0})\}$ and we obtain
\eqref{dplus+}. The other results follow similarly.
\end{proof}

Table \eqref{table1} shows the changes of $c_K$, $b_K$, $c_K - b_K + 1$ and
$c_K + b_K - 1$
under stabilization and destabilization of a closed braid.
\begin{equation}
\begin{tabular}{|l|l|l|l|l|}
\hline & $c_K$ & $b_K$ & $c_K-b_K+1$ & $c_K+b_K-1$ \\ \hline $+$ stabilization
& $+1$ & $+1$ & \multicolumn{1}{|c|}{$0$} &
\multicolumn{1}{|c|}{$+2$} \\ \hline $+$ destabilization & $-1$ &
$-1$ & \multicolumn{1}{|c|}{$0$} & \multicolumn{1}{|c|}{$-2$} \\
\hline $-$ stabilization & $-1$ & $+1$ &
\multicolumn{1}{|c|}{$-2$} & \multicolumn{1}{|c|}{$0$} \\ \hline
$-$ destabilization & $+1$ & $-1$ & \multicolumn{1}{|c|}{$+2$} &
\multicolumn{1}{|c|}{$0$} \\ \hline
\end{tabular}
\label{table1}
\end{equation}
Note that $c_K$ and $b_K$ are invariant under braid isotopy and exchange moves.

\begin{proof}[
Proof of \fullref{thmA}] Suppose that $K = K_\alpha
= K_+.$  Suppose we can apply positive destabilization $k$--times
($k\geq p$) to $K_{-}$. Let $\tilde{K}_{-}$ denote the closed
braid obtained after the destabilization. Then we have:

\begin{eqnarray}
d_{+}(P_{K_-}) + 2 &=& d_{+}(P_{\tilde{K}_-}) + 2  \notag \\
&\leq &( c_{\tilde{K}_-} + b_{\tilde{K}_-} -1 ) + 2  \notag \\
&=& \{(c_{K_-} + b_{K_-} -1) -2k \} + 2  \label{5eq} \\
&=&   (c_{K_+} -2) + b_{K_+} -1 -2k + 2   \notag \\
&=&   (c_{K_+} + b_{K_+} -1) -2k =  (c_K + b_K -1) -2k. \notag
\end{eqnarray}
The first equality holds since $K_{-}$ and $\tilde{K}_-$ have the
same knot type. The first inequality is the MFW--inequality. The
second equality follows from Table \eqref{table1}.

Similarly, if we can apply positive destabilization $l$--times $(l
\geq p)$ to $K_{0},$ and obtain $\tilde{K_0}$, we have
\begin{eqnarray}
d_{+}(P_{K_{0}})+1 &=& d_+ (P_{\tilde{K_0}}) + 1  \notag \\
&\leq& (c_{\tilde{K_0}} + b_{\tilde{K_0}} -1 ) + 1  \notag \\
&=&    (c_{K_0} + b_{K_0} -1 - 2l ) + 1  \label{5eq'} \\
&=&    (c_{K_+} -1 ) + b_{K_+} -1 - 2l + 1 \notag \\
&=&    (c_{K_+} + b_{K_+} -1) -2l =  (c_K + b_K -1) -2l. \notag
\end{eqnarray}
By \eqref{dplus+}, \eqref{5eq} and \eqref{5eq'} we get
\begin{eqnarray*}
d_+ (P_K) &=& d_+ (P_{K_+}) \leq \max \{d_{+}(P_{K_{-}})+2, \quad d_{+}(P_{K_{0}})+1\} \\
          &\leq& (c_K + b_{\cal K} -1) - \min \{ 2k, 2l \},
\end{eqnarray*}
ie, $D_K^+ \geq \min \{ 2k, 2l \} \geq 2p.$
When $K_\alpha = K_-$ or $K_\alpha = K_0$, the same arguments work
(use \eqref{dplus-} or \eqref{dplus0} for these cases in the place
of \eqref{dplus+}) and we get \eqref{less}.

The other inequality \eqref{more} also holds by the identical
argument.
\end{proof}

\subsection{Deficit growth}

Our goal is to exhibit examples (\fullref{deficit-thm})
of prime links on which the deficit
of the inequality can be arbitrary large.

\begin{figure}[ht!]\small
\labellist
\pinlabel {$K_+ =$} at 29 206
\pinlabel {$K_- =$} at 29 120
\pinlabel {$K_0 =$} at 29 35 
\pinlabel {negative} <-5pt,15pt> [lB] at 287 95
\pinlabel {destabilization} <-5pt,5pt> [lB] at 287 95
\pinlabel {positive} <-5pt,15pt> [lB] at 548 95
\pinlabel {destabilization} <-5pt,5pt> [lB] at 548 95
\pinlabel {positive} <-5pt,15pt> [lB] at 287 4
\pinlabel {destabilization} <-5pt,5pt> [lB] at 287 4
\endlabellist
\begin{center}\includegraphics [height=40mm]{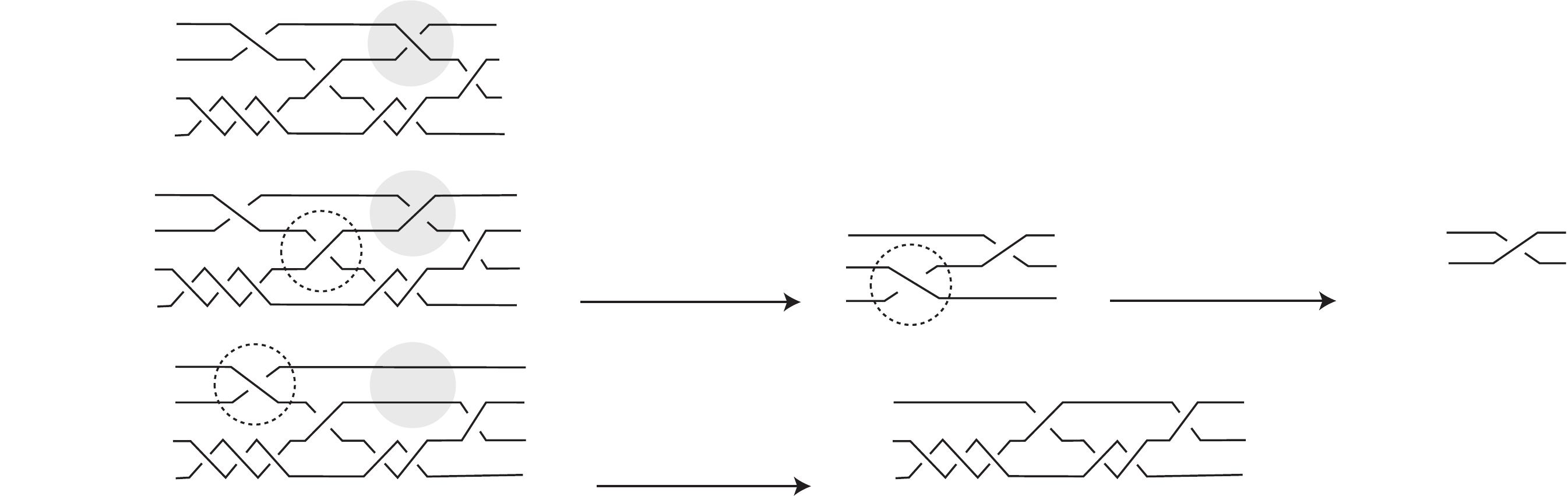}
\end{center}
\caption{Knot $9_{42}$ satisfies the conditions of \fullref{thmA}}\label{skein942}
\end{figure}
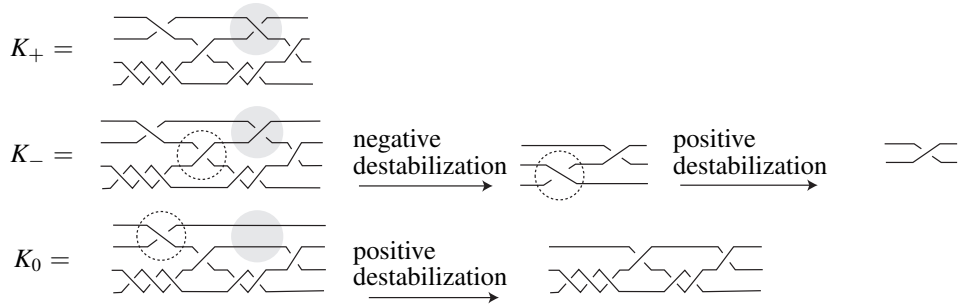

\begin{theorem}\label{th942}
Knot type ${\cal K}=9_{42}$ has a braid representative $K=K_+$
(see \fullref{skein942}) satisfying the sufficient
condition in \fullref{thmA}.
\end{theorem}

\begin{proof}[Proof of \fullref{th942}]
It is known that $9_{42}$ has braid index $=4$ and deficit $D_{9_{42}}=1.$
Let $K = K_+$ be its braid representative of the minimal braid index
as in \fullref{skein942}.
Construct $K_-, K_0$ by changing the shaded crossing.
Sketches show that both $K_-, K_0$ can be positively
destabilized. Thus by \fullref{thmA}, $D_{K}^+ \geq 2$ and
$D_{9_{42}}\geq 1.$
\end{proof}

\begin{theorem}\label{deficit-thm}
For any positive integer $n,$ there exists a prime link $L$ whose
deficit $D_L\geq n.$%
\end{theorem}
\begin{proof}[
Proof of \fullref{deficit-thm}] We prove the theorem
by exhibiting examples. For $n \in \mathbb{N}$ let ${\cal A}^n(9_{42})$ be
the closure of $n$--copies of $9_{42}$ linked each other by two
full twists as in the left sketch of \fullref{942}.
Since the braid index $b_{9_{42}} = 4$ and ${\cal A}^n(9_{42})$ is
an $n$--component link, we know the braid index of ${\cal A}^n(9_{42})$ is
$4n.$ This construction gives a braid representative with
$4n$--strands and $n$ distinguished crossings shaded in the left sketch.

\begin{figure}[ht!]
\begin{center}
\includegraphics [height=55mm]{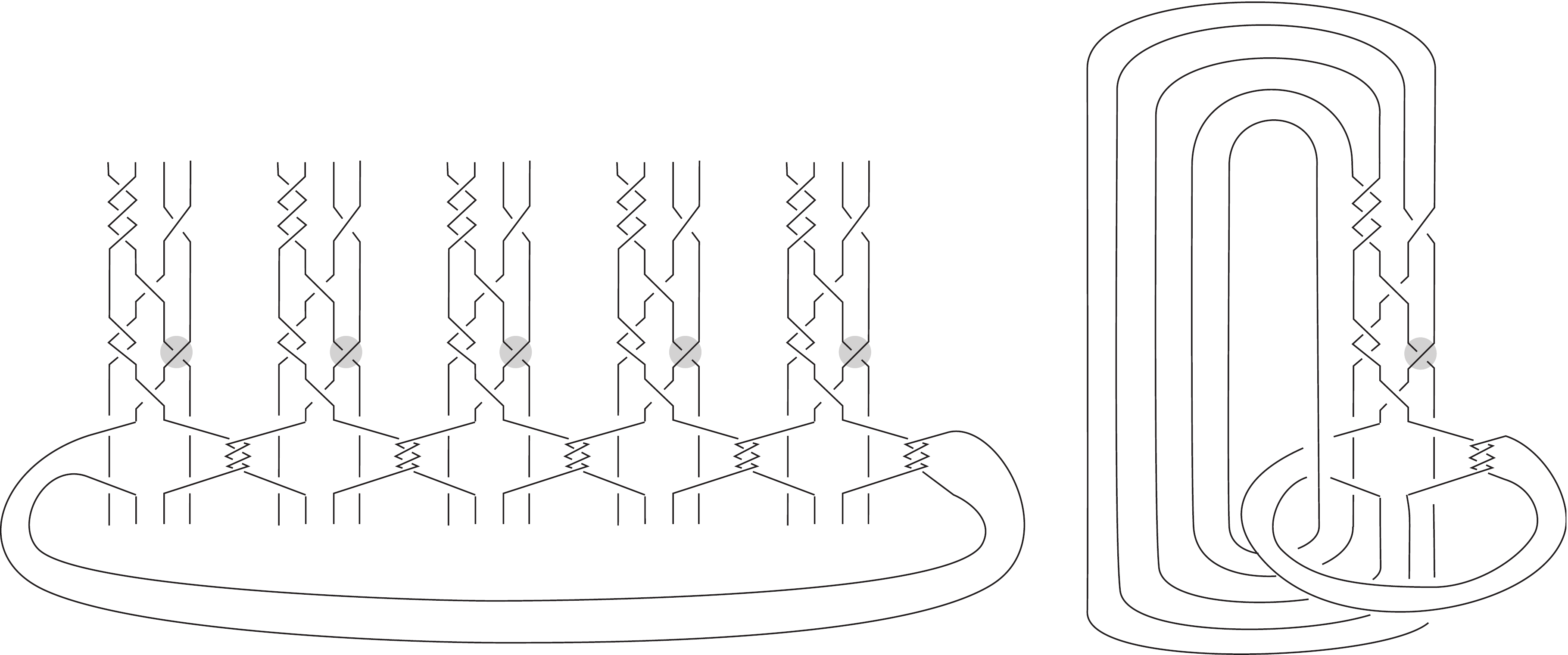}
\end{center}
\caption{Prime link ${\cal A}^5(9_{42})$
and $2$--component link ${\cal A}$} \label{942}%
\end{figure}

In the following we will see that each of the shaded crossing
contributes to the deficit.
\newpage

Let ${\cal K}:={\cal A}^2(9_{42})$ and let $K$ be the braid representative
of ${\cal K}$ as in \fullref{942}. Let $K_{- -}, K_{- 0}, K_{0
-}, K_{0 0}$ be the links obtained from $K$ by changing the two
shaded crossings. We repeat the discussion of the proof of \fullref{thmA}: We have:
\begin{eqnarray*}
d_+(P_{K_{- -}}) + (2+2) &=& d_{+}(P_{\tilde{K}_{- -}}) + 4   \\ %
&\leq &( c_{\tilde{K}_{- -}} + b_{\tilde{K}_{- -}} -1 ) + 4   \\
&=& \{(c_{K_{- -}} + b_{K_{- -}} -1) -2\cdot 2 \} + 4  \\
&=&   (c_{K} -4) + b_{K} -1 -2\cdot 2 + 4   \\
&=&   (c_{K} + b_{K} -1) -2\cdot 2
\end{eqnarray*}
Similarly:
\begin{eqnarray*}
 d_+(P_{K_{- 0}})+(2+1) &\leq & (c_{K} + b_{K} -1) -2\cdot 2 \\
 d_+(P_{K_{0 -}})+(1+2) &\leq & (c_{K} + b_{K} -1) -2\cdot 2 \\
 d_+(P_{K_{0 0}})+(1+1) &\leq & (c_{K} + b_{K} -1) -2\cdot 2
\end{eqnarray*}
Thus
\begin{eqnarray*}
d_+ (P_K) &=& \max \{d_+(P_{K_{- -}})+4, \  d_+(P_{K_{- 0}})+3, \ %
                     d_+(P_{K_{0 -}})+3, \  d_+(P_{K_{0 0}})+2 \} \\ %
          &\leq & (c_K + b_{\cal K} -1) - 2\cdot 2
\end{eqnarray*}
and\qquad\quad $D_{\cal K} \geq \frac{1}{2} D_K^+ \geq \frac{1}{2}(2\cdot 2) = 2.$

Similar arguments work when ${\cal K}={\cal A}^n(9_{42})$ for $n\geq 3$
and we have
 $D_{{\cal A}^n(9_{42})} \geq \frac{1}{2} D_{{\cal A}^n(9_{42})}^+
\geq \frac{1}{2} (2\cdot n)
 \geq n.$

The $2$--component link ${\cal A}$ of the right sketch is hyperbolic
\cite{N}. Pair $(S^3, {\cal A}^n(9_{42})\cup z\mbox{--axis})$ is an
$n$--fold cover of $(S^3, {\cal A})$ branched at $z$--axis . Therefore,
by Neumann and Zagier \cite{NZ} we can conclude that ${\cal
A}^n(9_{42})$'s are all prime except for finitely many $n$'s.
\end{proof}

\begin{remark}
{\rm By taking the connected sum of knots on which the MFW
inequality is non-sharp, one can also construct examples of
(non-prime) knots with arbitrarily large deficits. This fact
follows not only by \fullref{thmA} but also by
the definition of HOMFLYPT polynomial \eqref{skein1}
and the additivity of braid indices under connected sums
(Birman and Menasco \cite{BM4}).}
\end{remark}

\subsection{Birman--Menasco diagram}\label{section-BM}
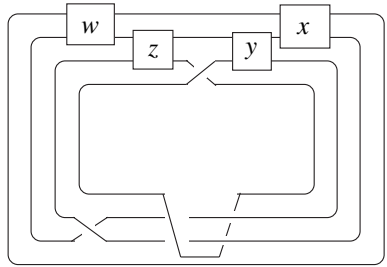
\begin{figure}[ht!]\small
\begin{center}
\labellist
\pinlabel {$w$} at 214 495
\pinlabel {$z$} at 271 474
\pinlabel {$y$} at 360 474
\pinlabel {$x$} at 407 495
\endlabellist
\includegraphics [height=35mm]{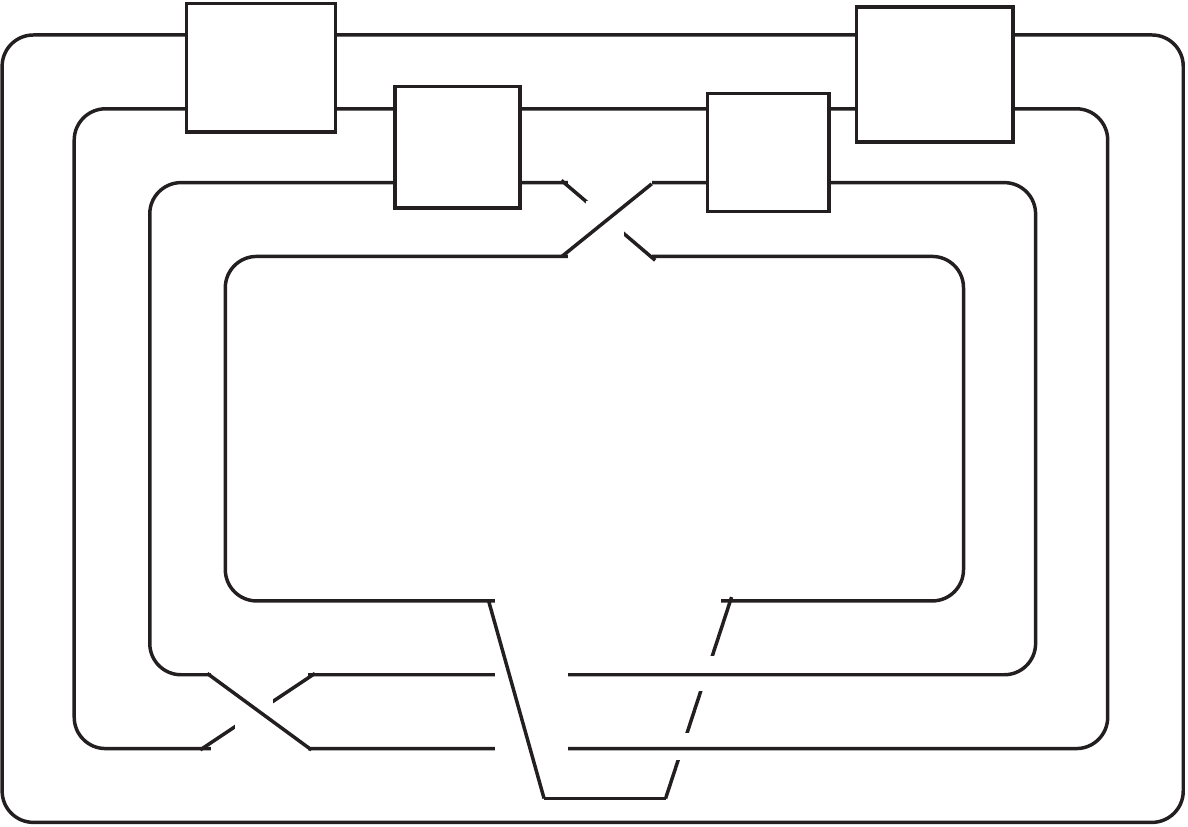}
\end{center}
\caption{The Birman--Menasco diagram $BM_{x,y,z,w}$} \label{menasco}
\end{figure}

As an application of \fullref{thmA}, we study
another infinite class of knots including all the Jones' five
knots  ($9_{42}, 9_{49}, 10_{132}, 10_{150}, 10_{156}$)
on which the MFW--inequality is not sharp.
We call the block-strand diagram (see \cite{MTWS-I} for definition)
of \fullref{menasco} the Birman--Menasco (BM) diagram.

\begin{definition}
Let $BM_{x,y,z,w}$, where $x,y,z,w \in \mathbb{Z}$, be the knot
(or the link) type which is obtained by assigning $x$ (resp.\ $y, z,
w)$ horizontal positive half-twists on two strands to the block $X$
(resp.\ $Y, Z, W)$ of the BM diagram.
\end{definition}

Recall that on all but only five knots ($9_{42}, 9_{49}, 10_{132},
10_{150}, 10_{156}$) up to crossing number $10$ the MFW--inequality
is sharp. An interesting property of the BM diagram is that it
carries all the five knots. Namely, we have
$9_{42}=BM_{-1,1,-2,-1}=BM_{-1,-2,-2,2}, \ $
$9_{49}=BM_{-1,1,1,2}\ $, $10_{132} = BM_{-1,-2,-2,-2}, \ $
$10_{150} = BM_{3,-2,-2,2} = BM_{-1,2,-2,2}$ $= BM_{-1,-2,2,2} =
BM_{-1,1,2,-1} =BM_{3,1,-2,-1},$ and $10_{156}=BM_{-1,1,1,-2}$.

We have the following theorem, which was conjectured informally by
Birman and Menasco: %
\begin{theorem}\label{BM-thm}
There are infinitely many $(x,y,z,w)$'s such that
the MFW--inequality is not sharp on $BM_{x,y,z,w}$.
\end{theorem}
\begin{lemma}\label{D+}
We have $D_{BM_{x,y,z,w}}^+ \geq 2$.
\end{lemma}

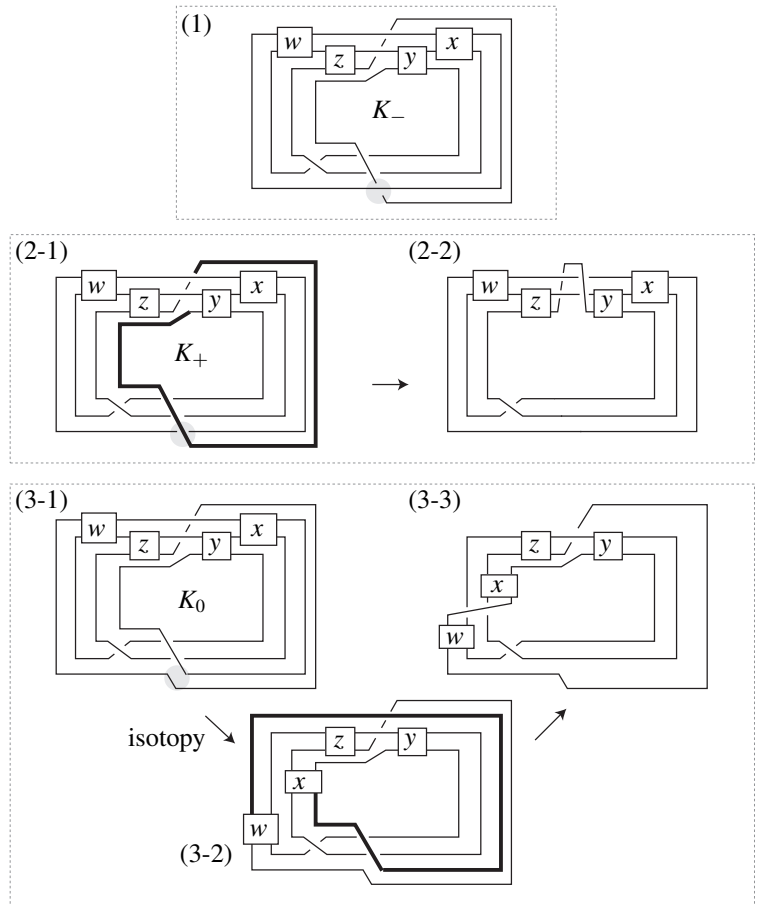
\begin{figure}[ht!]
\begin{center}
\labellist
\small
\pinlabel (1) [tl] at 130 703 
\pinlabel (2-1) [tl] at 0 525  
\pinlabel (2-2) [tl] at 308 525
\pinlabel (3-1) [tl] at 0 330
\pinlabel (3-2) [bl] at 129 27
\pinlabel (3-3) [tl] at 308 330
\pinlabel {$K_+$} at 143 431
\pinlabel {$K_-$} at 298 621
\pinlabel {$K_0$} at 143 239
\pinlabel {$x$} at 347 676
\pinlabel {$y$} at 314 662
\pinlabel {$z$} at 258 662
\pinlabel {$w$} at 222 676
\pinlabel {$w$} at 68 485
\pinlabel {$z$} at 105 472
\pinlabel {$y$} at 161 472
\pinlabel {$x$} at 194 485
\pinlabel {$w$} at 376 485
\pinlabel {$z$} at 412 472
\pinlabel {$y$} at 468 472
\pinlabel {$x$} at 500 485

\pinlabel {$w$} at 68 295
\pinlabel {$z$} at 105 282
\pinlabel {$y$} at 161 282
\pinlabel {$x$} at 194 295

\pinlabel {$w$} at 349 210
\pinlabel {$x$} at 382 251
\pinlabel {$z$} at 411 282
\pinlabel {$y$} at 467 281

\pinlabel {$w$} at 195 60
\pinlabel {$x$} at 227 97
\pinlabel {$z$} at 258 129
\pinlabel {$y$} at 314 129

\pinlabel {isotopy} [tr] at 158 142
\endlabellist 
\includegraphics [height=120mm]{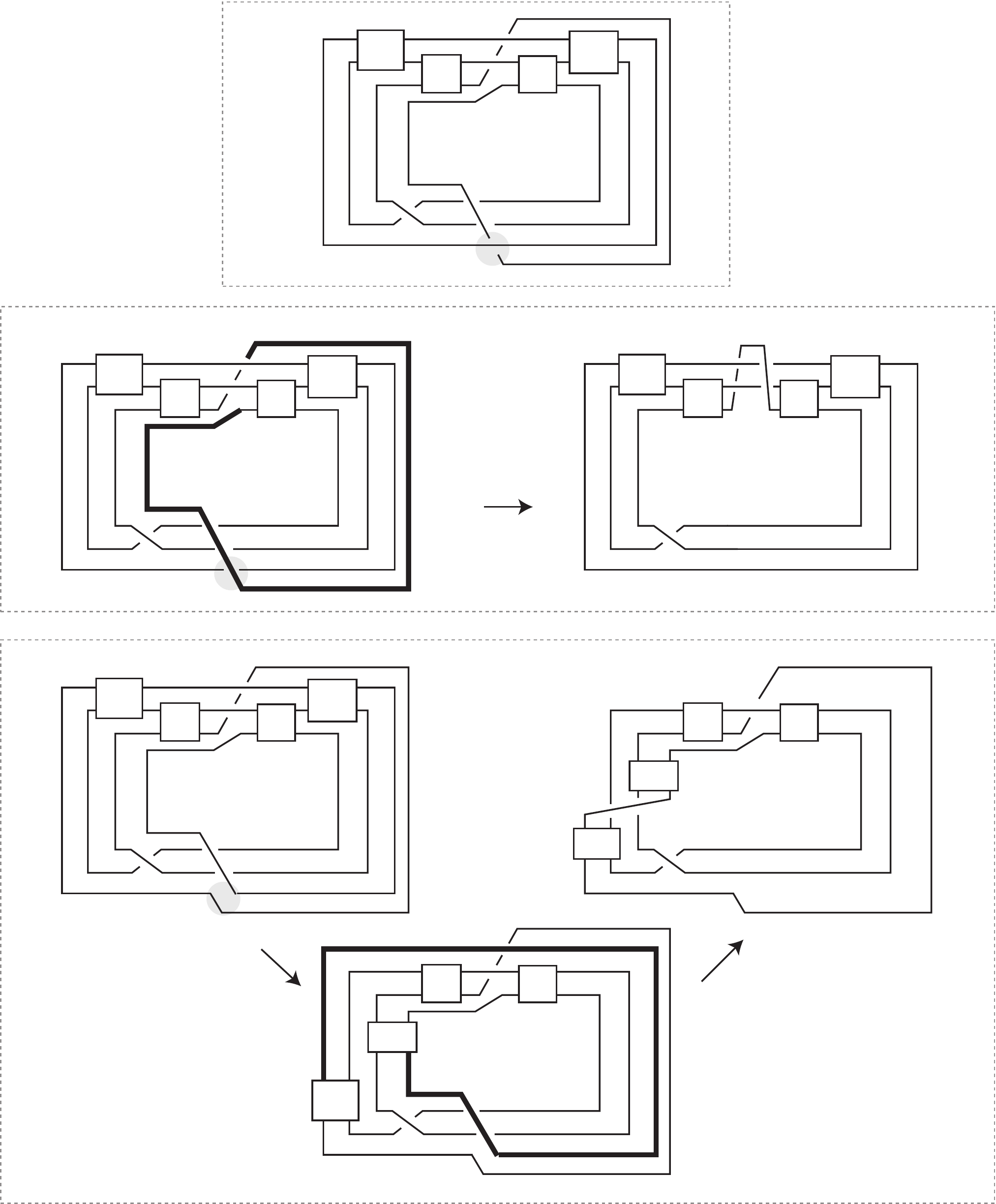}
\end{center}
\caption{The BM--diagram satisfies the sufficient condition} \label{menasco2}
\end{figure}

\begin{proof}[Proof of \fullref{D+}]
Change the BM diagram into the diagram in sketch (1) of \fullref{menasco2}%
by braid isotopy and denote it by $K$. Focus on the crossing
shaded in sketch (1).  Regard $K = K_-.$ We can apply positive
destabilization once to $K_+$ and obtain the diagram in sketch
(2-2). We also can apply positive destabilization once to $K_0$ as
we can see in the passage sketch (3-1) $\Rightarrow$
(3-2) $\Rightarrow$ (3-3). Therefore by \fullref{thmA}
we have $D_{BM_{x,y,z,w}}^+ \geq 2$ for any $(x,y,z,w)$.
\end{proof}

It remains to prove that there exist infinitely many
$(x,y,z,w)$'s such that the braid index of
$BM_{x,y,z,w}$ is $4$.  We introduce ${\cal K}_n :=
BM_{-1, -2, n, 2}$ and will show that for all $m \geq 1$ the
braid index of ${\cal K}_{2m}$ is $4$. (Note that ${\cal K}_{-2}=
9_{42}, {\cal K}_2 =10_{150}$ and ${\cal K}_{2m}$ is a knot.)
It will then follow, thanks to \fullref{D+}, that the MFW--inequality
cannot be sharp on any ${\cal K}_{2m},\ m\geq 1.$

In order to do this, we use the {\em enhanced Milnor number}
$\lambda$ defined by Neumann and Rudolph \cite{NR}.
Recall that the fiber surface of a fiber knot is obtained by
plumbing and deplumbing Hopf bands (see Giroux \cite{Giroux}).
This $\lambda$ is an invariant of fibered knots and links counting
algebraically the
number of negative Hopf bands to get the fiber surface.

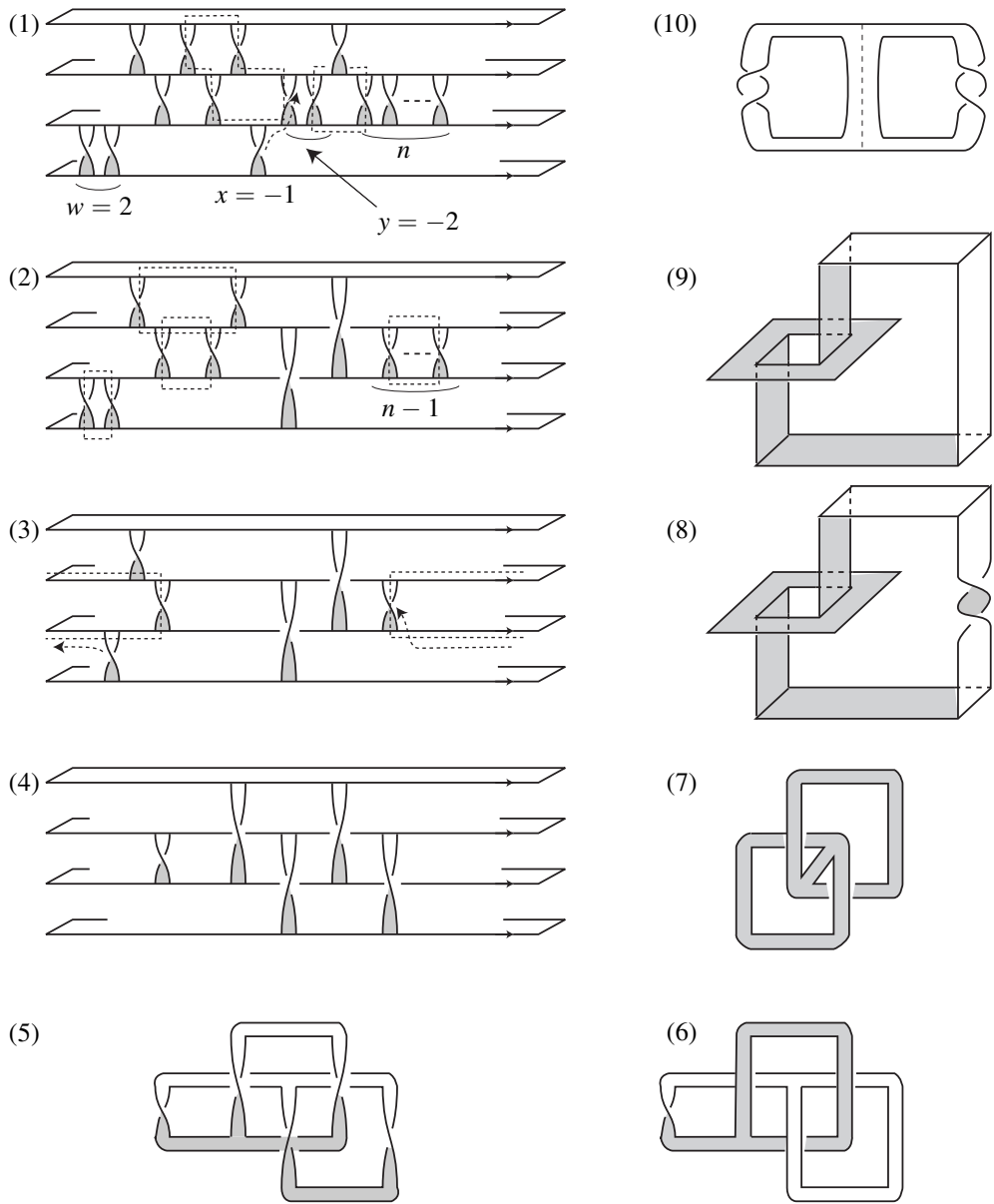
\begin{figure}[p]
\begin{center}
\labellist\small
\pinlabel (1) [r] at 35 760  
\pinlabel (2) [r] at 35 618  
\pinlabel (3) [r] at 35 476 
\pinlabel (4) [r] at 35 334 
\pinlabel (5) [r] at 35 194 
\pinlabel (10) [r] at 404 760  
\pinlabel (9) [r] at 404 618  
\pinlabel (8) [r] at 404 476 
\pinlabel (7) [r] at 404 334 
\pinlabel (6) [r] at 404 194 
\hair 1.5pt
\pinlabel {$x=-1$} [t] at 152 673  
\pinlabel {$y=-2$}  [tl] at 219 657 
\pinlabel {$w=2$} [t] at 64 666  
\pinlabel {$n$} [t] at 235 695
\pinlabel {$n-1$} [t] at 239 553
\endlabellist
\includegraphics [height=160mm]{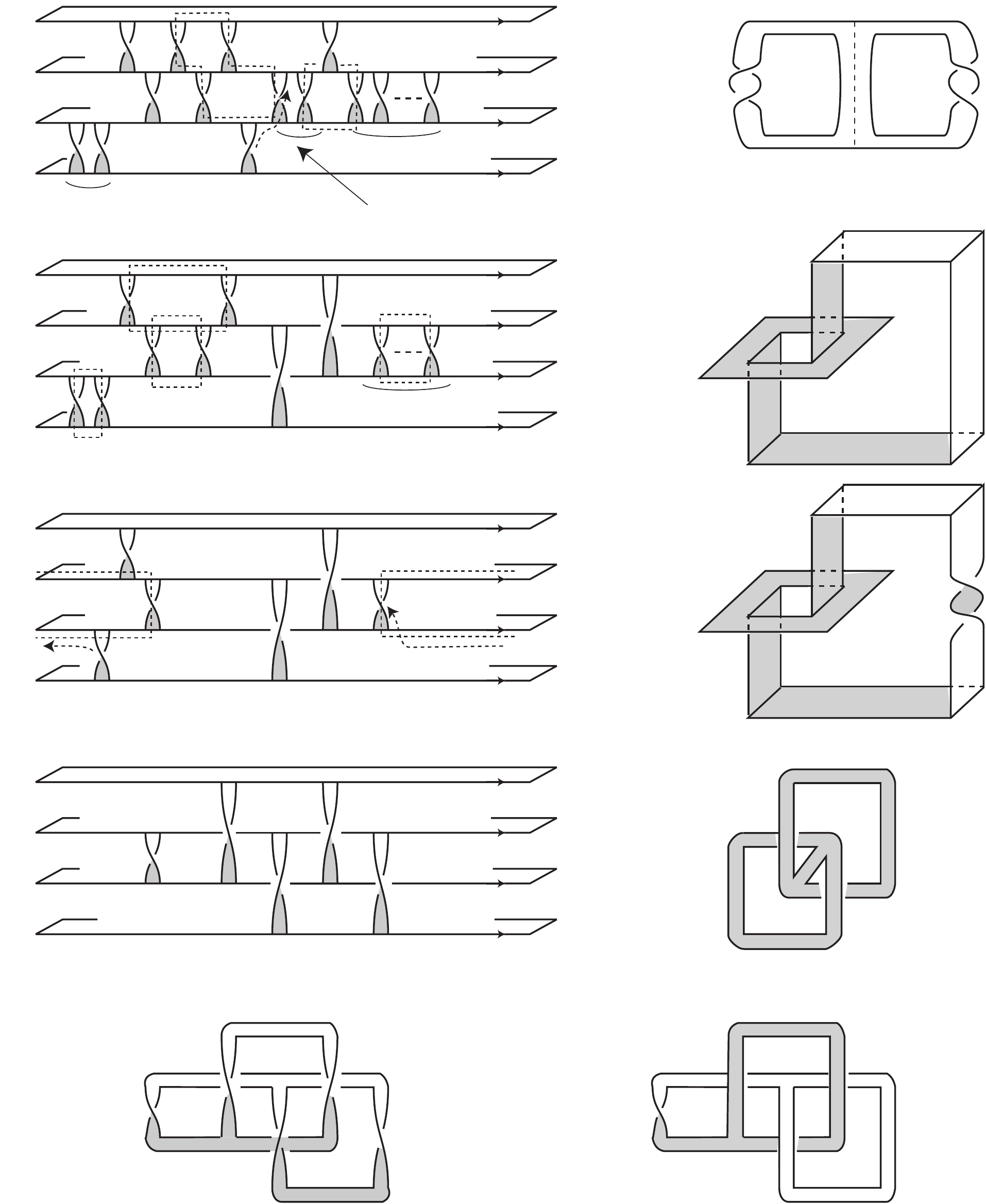}
\end{center}
\caption{Deformation of ${\cal K}_n$} \label{deplum}
\end{figure}
\newpage

\begin{lemma}\label{lambda=1}
All ${\cal K}_n$ $(n \geq 2)$ are fibered and have enhanced
Milnor number  $\lambda = 1.$
\end{lemma}
\begin{proof}[Proof of \fullref{lambda=1}]
Sketch (1) of \fullref{deplum} is the standard Bennequin surface of
${\cal K}_n$.
We compress it twice as in the passage sketch $(1) \Rightarrow (2)$
along the disks bounded by dotted circles in sketch (1).
Next, deplumb positive Hopf bands as much as possible as in the
passage sketch $(2) \Rightarrow (3) \Rightarrow (4)=(5)$.  Then
isotope the surface until we get sketch (8). These operations do
not change the enhanced Milnor number.

We apply Melvin and Morton's trick \cite{MM} p.167, as in the passage
sketch $(8) \Rightarrow (9).$
We remark that the enhanced Milnor number is invariant under this
trick.

The surface of sketch $(9)=(10),$ whose boundary is Pretzel link
$P(-2,0,2)$, is plumbing of a positive Hopf band and a negative
Hopf band. Thus it has $\lambda = 1$ so does ${\cal K}_n$.
\end{proof}

\begin{figure}[ht!]
\labellist\small
\pinlabel $a_1$ [t] at 121 460 
\pinlabel $a_2$ [t] at 324 460 
\pinlabel $a_3$ [t] at 492 460 
\endlabellist
\begin{center}
\includegraphics [height=24mm]{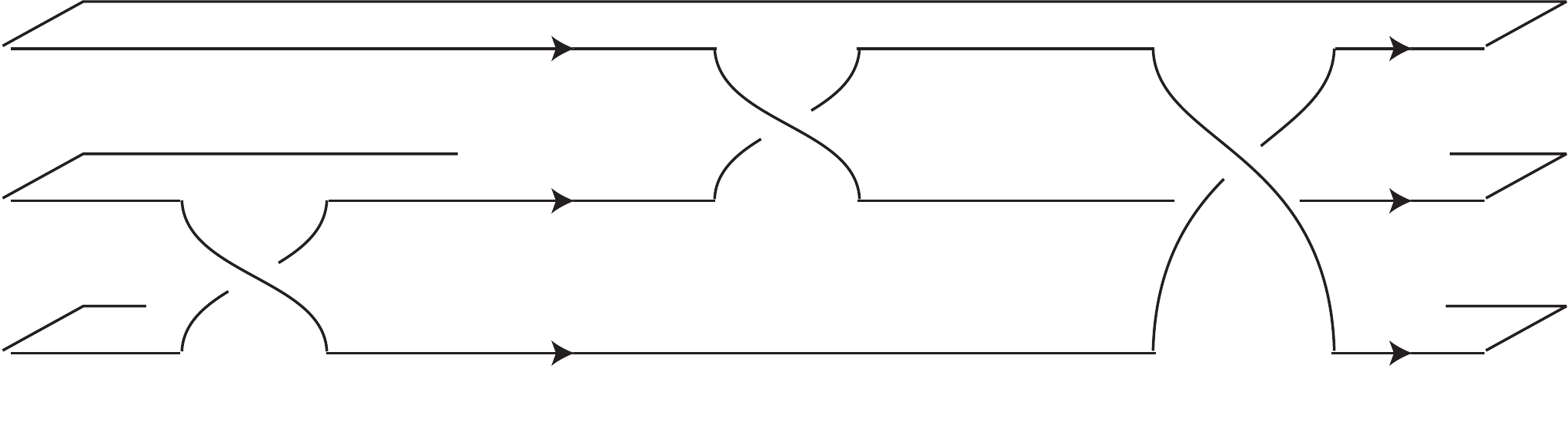}
\end{center}
\caption{Xu's band generators} \label{bands}%
\end{figure}

Here we summarize Xu's classification of $3$--braids \cite{Xu}. Let
$\sigma_1, \sigma_2$ be the standard generators of $B_3$ the braid
group of $3$--strings satisfying $\sigma_1 \sigma_2 \sigma_1 =
\sigma_2 \sigma_1 \sigma_2$. Let $a_1 := \sigma_1, a_2 :=
\sigma_2$ and $a_3 := \sigma_2 \sigma_1 \sigma_2^{-1}$. 
We can identify them with the twisted bands in \fullref{bands}.

Let
$\alpha := a_1 a_3 = a_2 a_1 = a_3 a_2$. If $w \in B_3$ let
$\overline{w}$ denote $w^{-1}.$

\begin{theorem}[Xu \cite{Xu}]\label{Xu-thm}
Every conjugacy  class in $B_3$ can be represented by a shortest
word in $a_1, a_2, a_3$ uniquely up to symmetry. And the word has
one of the three forms:
$$(1) \alpha^k P, \quad (2) N \overline{\alpha}^{k} , \quad (3) NP. $$
where $k\geq 0$ and $\overline{N}, P$ are positive words and the
arrays of subscripts of the words are non-decreasing.
\end{theorem}

The next is another lemma for \fullref{BM-thm}:
\begin{lemma}\label{ABCD}
If a closed $3$--braid has $\lambda=1$ and is a knot, then up to
symmetry it has one of the following forms:
\begin{eqnarray*}
A_x &:=& \overline{a_3} \ \overline{a_2} \ (a_1)^x, \quad x\geq 2, \mbox{\ even,} \\%
B_{x, y} &:=& \overline{a_3} \ \overline{a_3} \ (a_1)^x (a_2)^y,
                                     \quad x, y \geq 3, \mbox{\ odd,} \\%
C_{x, y, z} &:=& \overline{a_2} \ (a_1)^x (a_2)^y (a_3)^z,
                                     \quad x+z= \mbox{odd, } \ y= \mbox{even, } \ x,y,z\geq 1,\\ %
D_{x, y, z, w} &:=& \overline{a_2} \ (a_1)^x (a_2)^y (a_3)^z (a_1)^w, %
                         \quad x,y \geq 2, \ z,w\geq 1. %
\end{eqnarray*}
\end{lemma}
\begin{proof}[Proof of \fullref{ABCD}]
Assume we have a word $w \in B_3.$ By \fullref{Xu-thm}, $w$ has one
of the following forms:
\begin{center}
\begin{tabular}{|l|ll|}
\hline
%Case       & $w$             &             \\ \hline
\strutt Case (1)-1 & $w=\alpha^k$    & $k \geq 1$  \\ \hline
\strutt Case (1)-2 & $w=\alpha^k P$  & $k \geq 1$  \\ \hline
\strutt Case (1)-3 & $w=P$           & no $\alpha$ part \\ \hline
\strutt Case (2)-1 & $w=\overline{\alpha}^k$  & $k \geq 1$ \\ \hline
\strutt Case (2)-2 & $w=N\overline{\alpha}^k$ & $k \geq 1$ \\ \hline
\strutt Case (2)-2 & $w= N$          & no $\overline\alpha$ part \\ \hline
\strutt Case (3)   & $w=NP$          & $18$ cases to study \\ \hline
\end{tabular}
\end{center}

In this proof, we use the simplified notations:
\begin{center}
\begin{tabular}{|c|l|c|}
\hline
\strutt symbol & meaning & change in $\lambda$ \\ \hline
\strutt $i$    & $a_i$ for $i=1,2,3.$ & ---    \\ \hline
\strutt $=$    & same conjugacy class &  $0$   \\ \hline
\strutt $\longrightarrow$ & deplumb positive-Hopf bands & $0$  \\ \hline
\strutt $\Longrightarrow$ & deplumb negative-Hopf bands & $\geq 1$  \\ \hline
\strutt $\approx$         & Melvin--Morton trick \cite{MM} & $0$ \\ \hline
\strutt $\leadsto$ &composition of deplumbings of $\pm$ Hopf bands & $\geq 0$ \\ \hline
\end{tabular}
\end{center}

These are formulae we use:
\begin{eqnarray}
&& \alpha^2 \longrightarrow \alpha,  \mbox{ \ \fullref{alpha}} \label{1} \\
&& \alpha 1 2 3 \longrightarrow \alpha,  \mbox{ \ \fullref{alpha2}}
\label{2} \\
&& \overline{i} (i-1) i \ \approx \ \overline{i}\
\overline{i-1} i,   \mbox{ \ Melvin--Morton trick} \label{3} \\
&& i (i+1) \overline{i} \ \approx \ i \ \overline{i+1}\ \overline{i},
\mbox{ \ Melvin--Morton trick} \label{4}
\end{eqnarray}

\begin{figure}[ht!]
\begin{center}
\labellist\small\hair 1.5pt
\pinlabel $\alpha$ [t] at 61 695
\pinlabel $\alpha$ [t] at 107 695
\pinlabel $\alpha$ [t] at 541 695
\hair 0pt
\pinlabel deplumb [l] <0pt,-10pt> at 146 685
\pinlabel slide [l] <0pt,-10pt> at 317 685
\pinlabel deplumb [l] <0pt,-10pt> at 459 685
\endlabellist 
\includegraphics [width=.98\hsize]{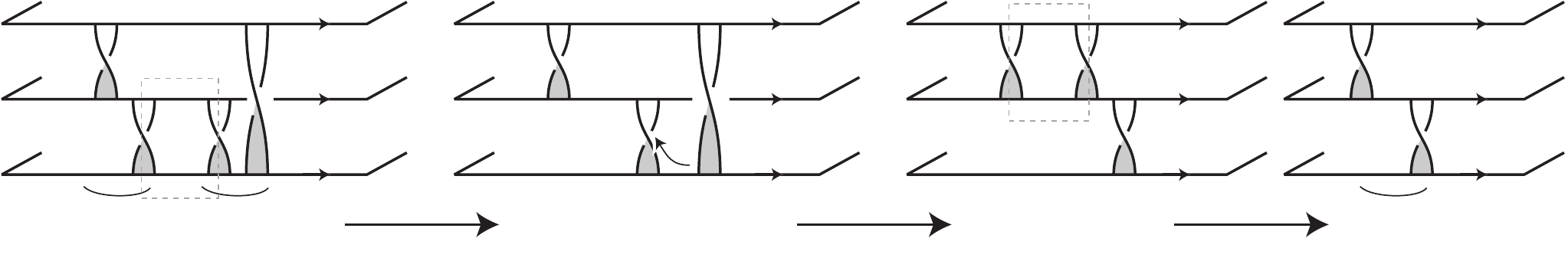}
\end{center}
\caption{$\alpha^2 \longrightarrow \alpha$} \label{alpha}%
\end{figure}

\begin{figure}[ht!]
\begin{center}
\labellist\small\hair 1.5pt
\pinlabel $\alpha$ [t] at 67 366
\pinlabel $\alpha$ [t] at 515 366
\hair 0pt
\pinlabel deplumb [l] <0pt,-10pt> at 145 356
\pinlabel deplumb [l] <0pt,-10pt> at 285 356
\pinlabel deplumb [l] <0pt,-10pt> at 432 356
\pinlabel slide [l] <0pt,-20pt> at 145 356
\pinlabel slide [l] <0pt,-20pt> at 285 356
\endlabellist 
\includegraphics [width=.98\hsize]{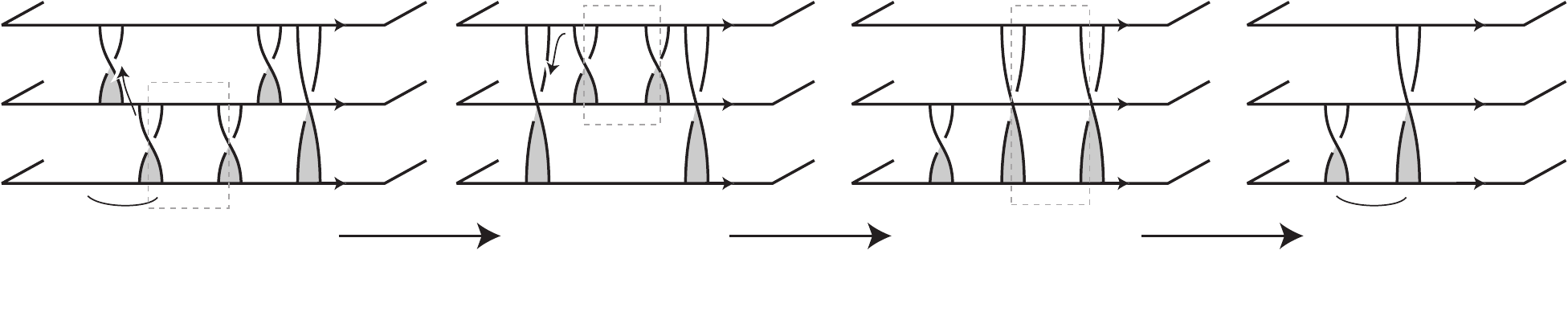}
\end{center}
\caption{$\alpha 1 2 3 \longrightarrow \alpha$}
\label{alpha2}\end{figure}

Now we study each case.

\noindent\textbf{Case (1)-1}\qua
By \eqref{1},
we have $w=\alpha^k \longrightarrow \alpha (=$ unknot).
Thus $w$ has $\lambda=0.$

\noindent\textbf{Case (1)-2}\qua By \eqref{1} up to permutation of
$\{1,2,3\}$ we have $\alpha^k P \longrightarrow \alpha
P \longrightarrow$\break $\alpha (1 2 3 1 2 3 \cdots \cdots).$ Thanks to
\eqref{2} we have
$$\alpha \ \overbrace{1 2
3 1 2 3 \cdots \cdots}^{\mbox{length$=l$}} \ \longrightarrow \
\alpha \ \overbrace{1 2 3 1 2 3 \cdots
\cdots}^{\mbox{length $=l-3$}} \quad \mbox{for } l\geq 3.$$ If
$l=1, 2$, we have $\alpha 1 = 2 1 1 \longrightarrow \alpha$ and
$\alpha 1 2 = 2 1 1 2 \longrightarrow \alpha$. Thus $w$ has
$\lambda=0.$

\noindent\textbf{Case (1)-3}\qua Assume $w=P$. There are three possible cases:
$$P \longrightarrow (1 2 3)^n, \quad
P \longrightarrow (1 2 3)^n 1 \quad \mbox{ and } \quad P
\longrightarrow (1 2 3)^n 1 2 \quad \mbox{where}\ n\geq 0.$$

If $w$ satisfies the first case, it is proved that $(123)^n$ is not fibered
in Theorem 3.2 of \cite{S-preprint}, where Stoimenow determines fibreness of 
strongly quasi-positive $3$--braid links. Therefore, $w$ is not fibered.

The second case can be reduced to the first case, since $(1 2 3)^n
1 = 1 (1 2 3)^n \longrightarrow (1 2 3)^n.$

For the third case, since $(1 2 3)^n 1 2 =  2 (1 2 3)^n 1 = \alpha
(2 3 1)^n \longrightarrow \alpha,$ $w$ has $\lambda=0.$

\noindent\textbf{Case (2)-1}\qua
By \eqref{1}, $w=\overline{\alpha}^k\Longrightarrow\overline{\alpha}$
and $w$ has $\lambda=2(k-1)\neq 1$.

\noindent\textbf{Case (2)-2}\qua
Suppose $w=N \overline{\alpha}^k$ where $k\geq 1$.

If $w=\overline{i}\ \overline{\alpha},$ we have $\overline{i}\
\overline{\alpha} \Longrightarrow \overline{\alpha}$ and $w$ has
$\lambda=1.$ However, the closure of $w$ has more than one
component and it does not satisfy the condition of the lemma.

If $w\neq \overline{i}\ \overline{\alpha},$ we have $N
\overline{\alpha}^k \Longrightarrow \overline{\alpha}$ by
\eqref{2}, and $w$ has $\lambda \geq 2.$

\noindent\textbf{Case (2)-3}\qua
Suppose $w=N$. There are three possible cases:
$$N \Longrightarrow (\overline{3} \ \overline{2}\
\overline{1})^n, \quad N \Longrightarrow (\overline{3} \
\overline{2}\ \overline{1})^n \overline{3} \quad \mbox{and} \quad
N \Longrightarrow (\overline{3} \ \overline{2}\ \overline{1})^n
\overline{3} \ \overline{2}\quad \mbox{where}\ n\geq 0.$$

For the first case, $w$ is not fibered \cite{S-preprint}.

For the second case, if $n=0$ then  $w$ has $\lambda=1$ if and
only if $w=\overline{3}\ \overline{3}$. However this has two
components. If $n\geq 1,$ since $(\overline{3} \ \overline{2}\
\overline{1})^n \overline{3} \Longrightarrow (\overline{3} \
\overline{2}\ \overline{1})^n$ it can be reduced to the first
case.

For the third case, if $n=0$ then $w$ has $\lambda=1$ if and only
if $w=\overline{3}\ \overline{3} \ \overline{2}.$ However it has
two components. If $n\geq 1,$ we have $(\overline{3} \
\overline{2} \ \overline{1})^n
  \overline{3} \ \overline{2}
= \overline{2} \ \overline{3}
  (\overline{2} \ \overline{1} \ \overline{3})^n
=\overline{\alpha}\ (\overline{2} \ \overline{1} \ \overline{3})^n
\Longrightarrow \overline{\alpha}$ and $w$ has $\lambda \geq 3n.$

\noindent\textbf{Case (3)}\qua Assume $w=NP.$
Let $w'$ be a word obtained from $w$
by deplumbing $\pm$ Hopf bands sufficiently enough times, ie, $w \leadsto
w'.$ This $w'$ has one
of the following 18 forms up to permutation of $\{1,2,3\}$.

$$
\begin{tabular}{|l|lll|}
\hline case & $w'$ & & \\ \hline
\strutt i & $(\overline{2}\ \overline{1}\ \overline{3}%
)^{k}(123)^{l}$ & $k, l \geq 1$ & \\ \hline
\strutt ii & $(\overline{2}\ \overline{1}\ \overline{3}%
)^{k}(123)^{l}1$ & $k\geq 1,l\geq 0$ & \\ \hline
\strutt iii & $(\overline{2}\ \overline{1}\ \overline{3}%
)^{k}(123)^{l}12$ & $k\geq 1,l\geq 0$ & not shortest word \\ \hline
\strutt iv & $\overline{3}(\overline{2}\ \overline{1}\ \overline{3}%
)^{k}(123)^{l}$ & $k\geq 0,l\geq 1$ & not shortest word \\ \hline
\strutt v & $\overline{3}(\overline{2}\ \overline{1}\ \overline{3}%
)^{k}(123)^{l}1$ & $k, l \geq 0$ & \\ \hline
\strutt vi & $\overline{3}(\overline{2}\ \overline{1}\ \overline{3}%
)^{k}(123)^{l}12$ & $k, l \geq 0$ & \\ \hline
\strutt vii & $\overline{1}\ \overline{3}(\overline{2}\ \overline{1}%
\ \overline{3})^{k}(123)^{l}$ & $k\geq 0,l\geq 1$ &
\\ \hline
\strutt viii & $\overline{1}\ \overline{3}(\overline{2}\ \overline{1}%
\ \overline{3})^{k}(123)^{l}1$ & $k, l \geq 0$ & not shortest word \\
\hline \strutt ix & $\overline{1}\ \overline{3}(\overline{2}\
\overline{1}\
\overline{3})^{k}(123)^{l}12$ & $k, l\geq 0$ & \\
\hline
\end{tabular}$$
$$
\begin{tabular}{|l|lll|}
%\hline case & $w'$ & &  \\
\hline
\strutt i$'$ & $(\overline{1}\ \overline{3}\ \overline{2}%
)^{k}(123)^{l}$ & $k, l\geq 1$ & \\ \hline
\strutt ii$'$ & $(\overline{1}\ \overline{3}\ \overline{2}%
)^{k}(123)^{l}1$ & $k\geq 1,l\geq 0$  & not shortest word \\ \hline
\strutt iii$'$ & $(\overline{1}\ \overline{3}\ \overline{2}%
)^{k}(123)^{l}12$ & $k\geq 1,l\geq 0$ & \\
\hline
\strutt iv$'$ & $\overline{2}(\overline{1}\ \overline{3}\ \overline{2}%
)^{k}(123)^{l}$ & $k\geq 0,l\geq 1$ & \\ \hline
\strutt v$'$ & $\overline{2}(\overline{1}\ \overline{3}\ \overline{2}%
)^{k}(123)^{l}1$ & $k, l\geq 0$ & \\ \hline
\strutt vi$'$ & $\overline{2}(\overline{1}\ \overline{3}\ \overline{2}%
)^{k}(123)^{l}12$ & $k, l\geq 0$  & not shortest word \\ \hline
\strutt vii$'$ & $\overline{3}\ \overline{2}(\overline{1}\ \overline{3}%
\ \overline{2})^{k}(123)^{l}$ & $k\geq 0,l\geq 1$
 & not shortest word \\ \hline
\strutt viii$'$ & $\overline{3}\ \overline{2}(\overline{1}\ \overline{3%
}\ \overline{2})^{k}(123)^{l}1$ & $k, l\geq 0$ & \\ \hline
\strutt ix$'$ & $\overline{3}\ \overline{2}(\overline{1}\ \overline{3}%
\ \overline{2})^{k}(123)^{l}12$ & $k, l\geq 0$ & \\ \hline
\end{tabular}$$
Since words of case iii, iv, viii, ii$'$, vi$'$ and vii$'$ are not shortest
(reducible) we eliminate them from the list.

These are reduction formulae we use:
\begin{eqnarray}
  (\ol{2}\ \ol{1}\ \ol{3})(123) &\stackrel{\eqref{F}}{\lra} & \ol{2}(3\ol{2})
  = 1\ol{2}\ \ol{2} \Lra 1\ol{2}
\label{reduction1} \\
(1 \ol{2})(123) &\stackrel{\eqref{S}}{\lra}& 11 \ol{2}\lra 1 \ol{2}
\label{A}\\
(\overline{2}\ \overline{1}\ \overline{3})(1 \ol{2})
&\stackrel{\eqref{C}}{\lra} & 1\ol{3}\ \ol{2} = \ol{2}\ \ol{2}3 \Lra\ =1 \ol{2}
\label{B} \\
(\ol{2}\ \ol{1}\ \ol{3})1 &\approx & \ol{2}\ \ol{1}31=1\ol{2}\ \ol{2}1\Lra\ =
11\ol{3}\lra  1 \ol{3}
\label{C}  \\
(\ol{2}\ \ol{1}\ \ol{3})1 \ol{3} &\stackrel{\eqref{C}}{\lra}&  1 \ol{3}\ \ol{3}
\Lra 1 \ol{3}
\label{D} \\
\ol{1}\ \ol{3}(123) &=& \ol{1}22\ol{1}3\lra \ol{1}2\ol{1}3 =3\ol{1}\ \ol{1}3
\Lra\ = 33\ol{2} \lra 3 \ol{2}
\label{F} \\
3\ol{2}(123)  & \stackrel{\eqref{S}}{\lra}& 3(1\ol{2}) = \ol{1}33\lra \ol{1}3
= 3\ol{2}
\label{G}  \\
(123)1\ol{3} &\approx & 123\ol{1}\ \ol{3} = 1\ol{3}22\ol{3} \nonumber \\
&\lra & 1\ol{3}2\ol{3} = 1\ol{3}\ \ol{3}1\Lra 1\ol{3}1\lra 1 \ol{3}
\label{H} \\
12 \ol{1}\ \ol{3} &=& 1\ol{3}\ \ol{3}1\Lra 1\ol{3}1 \lra \ol{2}1=1 \ol{3}
\label{N} \\
(\ol{1}\ \ol{3}\ \ol{2})(123) &\stackrel{\eqref{Q}}{\lra}& 3\ol{1}3=33\ol{2}
\Lra 3\ol{2}=\ol{1}3
\label{reduction2} \\
(\ol{1}\ \ol{3}\ \ol{2}) \ol{1}3 &=& \ol{1}\ \ol{3}\ 1 \ol{2}\ \ol{2}
\Lra\ \approx \ol{1} 31 \ol{2} = \ol{1}\ \ol{1} 33 \lra\ \Lra \ol{1} 3
\label{O} \\
\ol{1}3 (123) &\approx & \ol{1}\ \ol{3} 123 = \ol{1}22\ol{1}3  \lra\ \Lra
3\ol{1}3 \lra 3 \ol{2} = \ol{1}3
\label{P} \\
(\ol{1}\ \ol{3}\ \ol{2})12 &=& \ol{1}2\ol{3}\ \ol{3}2\Lra\ =  \ol{1}22\ol{1}
\lra \ol{1}2\ol{1}=3\ol{1}\ \ol{1} \Lra  3\ol{1}
\label{Q} \\
(\ol{1}\ \ol{3}\ \ol{2})3\ol{1} &\stackrel{\eqref{B}}{\lra} & 3\ol{1} \quad
\mbox{permutation of \eqref{B}}
\label{R} \\
\ol{2}(123) &\approx & \ol{2}\ \ol{1}23  = \ol{2} 33\ol{2} \lra \ol{2}
3\ol{2} = 1 \ol{2} \ \ol{2} \Lra 1 \ol{2}
\label{S} \\
1 \ol{2}\ (\ol{1}\ \ol{3}\ \ol{2}) &\approx & 12\ol{1}\ \ol{3}\ \ol{2}=
1\ol{3}\ \ol{3}1\ol{2} \Lra\ = \ol{2}11\ol{2} \lra  \ol{2}1\ol{2}
\Lra 1 \ol{2}
\label{W}
\end{eqnarray}

\begin{sublemma}\label{red3}
For $k, l \geq 1,$ we have:
\begin{eqnarray}
(\overline{2}\ \overline{1}\ \overline{3})^{k}
(123)^l & \leadsto & 1 \overline{2} \label{AA} \\
(\ol{1}\ \ol{3}\ \ol{2})^k (123)^l & \leadsto & \ol{1}3 \label{BB}
\end{eqnarray}
Either case, the increase of $\lambda$ is $\geq 2.$
\end{sublemma}

\begin{proof}
From \eqref{reduction1}, \eqref{B} and \eqref{C}, 
we obtain \eqref{AA}. Similarly, \eqref{BB} follows from
\eqref{reduction2}, \eqref{O} and \eqref{P}.
\end{proof}

\noindent\textbf{Case 3-i}\qua
By \eqref{AA}, our $w$ is fibered and $\lambda \geq 2.$

\noindent\textbf{Case 3-ii}\qua
If $k\geq 1. l=0,$
$$w \leadsto w' = (\overline{2}\ \overline{1}\ \overline{3})^{k} 1
\ \stackrel{\eqref{C}}{\lra} \
(\overline{2}\ \overline{1}\ \overline{3})^{k-1} 1 \ol{3} \
\stackrel{\eqref{D}}{\lra} \cdots \stackrel{\eqref{D}}{\lra} \ 1 \ol{3}
$$
and $w$ has $\lambda = 1$ if and only if $w=\ol{2}\ \ol{1}\
\ol{3}\ 1^x$ for some $x \geq 1.$ However it has $2$ or
$3$ components and it does not satisfy the condition of \fullref{ABCD}.
If $k, l \geq 1,$
$$w \leadsto w' = (\overline{2}\ \overline{1}\ \overline{3})^{k}
(123)^l 1 \ \stackrel{\eqref{AA}}{\lra} \ (1 \ol{2}) 1 \lra 1 \ol{2}$$
and $w$ has $\lambda \geq 2.$

\noindent\textbf{Case 3-v}\qua
When $k=l=0,$ $w$ has $\lambda=1$ if and only if
$w= \ol{3}\ \ol{3}\ 1^x$ for some $x\geq 1,$ which has more than $1$ component.
If $k \geq 1$ and $l=0$,
$$w \leadsto w' = \ol{3}\ (\ol{2}\ \ol{1}\ \ol{3})^k 1
\ \stackrel{\eqref{D}}{\lra}\ \cdots\ \stackrel{\eqref{D}}{\lra}\ 1 \ol{3}$$
and $\lambda \geq 2.$
If $k=0, l\geq 1,$
$$w\leadsto w'=\ol{3}\ (123)^l1=(123)^l1\ol{3}\ \stackrel{\eqref{H}}{\lra}\
\cdots\ \stackrel{\eqref{H}}{\lra} \ 1 \ol{3}.$$
Thus $w$ has $\lambda=1$ if and only if $w= \ol{3}\  1^x\ 2^y\ 3^z\ 1^w$
for $x,y,z,w \geq 1.$
If $k, l \geq 1,$
$$w \leadsto w' = \ol{3}\ (\ol{2}\ \ol{1}\ \ol{3})^k (123)^l 1
\ \stackrel{\eqref{AA}}{\lra} \ \ol{3} (1 \ol{2}) 1 =
\ol{3} 11 \ol{3}\lra\ \Lra 1 \ol{3}$$ and $\lambda \geq 2.$

\noindent\textbf{Case 3-vi}\qua
If $k=l=0$,
$$w \leadsto w' = \ol{3}\ 1 2 = 22\ol{1} \lra 2 \ol{1}.$$
Therefore, $w$ has $\lambda=1$ if and only if $w= \ol{3}\ \ol{3}\ 1^x\ 2^y$
for $x,y \geq 1.$
If $k=0, l\geq 1$
$$w \leadsto w' = \ol{3} (123)^l 1 2 =(123)^l12\ol{3}\lra (123)^l 1 \ol{3}\
\stackrel{\eqref{H}}{\lra} \ \cdots \ \stackrel{\eqref{H}}{\lra} \ 1 \ol{3}$$
and $w$ has $\lambda=1$ if and only if $w= \ol{3} 1^x\ 2^y\ 3^z\ 1^w\ 2^v$
for some $x,y,z,w,v \geq 1.$
When $k\geq 1, l=0$,
$$w \leadsto w' = \ol{3} (\ol{2}\ \ol{1}\ \ol{3})^k 12 \lra (\ol{2}\ \ol{1}\
\ol{3})^k 1 \ol{3}  \stackrel{\eqref{D}}{\lra} \ \cdots \
\stackrel{\eqref{D}}{\lra} \ 1 \ol{3}$$ and $\lambda \geq 2$.
If $k, l \geq 1,$
$$w \leadsto w'\ \stackrel{\eqref{AA}}{\lra}\ \ol{3}(1\ol{2})12 =\ol{3}11\ol{3}
2 \lra\ol{3}1\ol{3} 2 = 2 \ol{3}\ \ol{3}\ 2 \lra\ \Lra 2 \ol{3}$$ and
$\lambda \geq 2$.

\noindent\textbf{Case 3-vii}\qua
If $k=0, l \geq 1,$
$$w \leadsto w' = \ol{1}\ \ol{3}\ (123)^l\ \stackrel{\eqref{F}}{\lra} \
3 \ol{2} (123)^{l-1} \ \stackrel{\eqref{G}}{\lra} \ \cdots
\ \stackrel{\eqref{F}}{\lra} \ 3 \ol{2}$$ and $\lambda =1$ if and
only if $w=\ol{1}\ \ol{3}\ 1^x\ 2^y\ 3^z$ for some $x,y,z \geq 1.$
 If $k, l \geq 1,$  $$w \leadsto w'\ \stackrel{\eqref{AA}}{\lra} \
\ol{1}\ \ol{3}(1 \ol{2})=\ol{2}\ \ol{1}\ \ol{3}1\ \stackrel{\eqref{C}}{\lra}\
1 \ol{3}$$ and $\lambda \geq 2$.

\noindent\textbf{Case 3-ix}\qua
When $k=l=0,$
$$w \leadsto w' = \ol{1}\ \ol{3}\ 12 = \ol{1}22\ol{1}
\lra \ \Lra  \ol{1}2.$$ Thus $w$ has $\lambda=1$ if and only if
$w= \ol{1}\ \ol{3}\ 1^x\ 2^y$ for some $x,y \geq 1.$
When $k=0, l\geq 1$
$$w \leadsto w' =  \ol{1}\ \ol{3} (123)^l 12
\ \stackrel{\eqref{F}}{\lra}  \ \stackrel{\eqref{G}}{\lra} \
\ol{1} 312 = 3 \ol{2} 12 \approx 3 \ol{2}\ \ol{1} 2 = \ol{1} 33 \ol{1}
\lra\ \Lra \ol{1}3$$ and $\lambda\geq 2.$
When $k\geq 1, l=0$
$$w \leadsto w' =  \ol{1}\ \ol{3} (\ol{2}\ \ol{1}\ \ol{3})^k 12
= 12 \ol{1}\ \ol{3} (\ol{2}\ \ol{1}\ \ol{3})^k  \ \stackrel{\eqref{N}}{\lra} \
1 \ol{3} (\ol{2}\ \ol{1}\ \ol{3})^k \ \stackrel{\eqref{D}}{\lra} \ \cdots
 \ \stackrel{\eqref{D}}{\lra} \ 1 \ol{3}$$ and  $\lambda\geq 2.$
When $k,l\geq 1$
$$w \leadsto w' \ \stackrel{\eqref{AA}}{\lra} \
\ol{1}\ \ol{3} (1\ol{2}) 12 = 12\ol{1}\ \ol{3}1\ol{2}
\ \stackrel{\eqref{N}}{\lra} \ 1 \ol{3}1\ol{2}  = \ol{2} 11 \ol{2}
\lra\ \Lra 1 \ol{2}$$  and  $\lambda\geq 2.$

\noindent\textbf{Case 3-i$'$}\qua
By \eqref{BB} our $w$ is fibered and $\lambda \geq 2.$

\noindent\textbf{Case 3-iii$'$}\qua
When $k>l=0,$
$$w \leadsto w' = (\ol{1}\ \ol{3}\ \ol{2})^k 12 \ \stackrel{\eqref{Q}}{\lra} \
(\ol{1}\ \ol{3}\ \ol{2})^{k-1} 2 \ol{1}   \ \stackrel{\eqref{R}}{\lra} \ \cdots
 \ \stackrel{\eqref{R}}{\lra} \ 3 \ol{1}$$
and $\lambda \geq 2.$
When $k,l\geq 1$,
$$w \leadsto w' \ \stackrel{\eqref{BB}}{\lra} \ (\ol{1}3)12 \approx
\ol{1}\ \ol{3}\ 12 = \ol{1}22\ol{1} \lra\ \Lra 3 \ol{1}$$
\newpage
and $\lambda \geq 2$.

\noindent\textbf{Case 3-iv$'$}\qua
When $k=0, l\geq 1,$
$$w \leadsto w'= \ol{2}(123)^l  \ \stackrel{\eqref{S}}{\lra} \
1 \ol{2}(123)^{l-1} \ \stackrel{\eqref{A}}{\lra} \ \cdots
\ \stackrel{\eqref{A}}{\lra} \ 1 \ol{2}.$$
Thus $w$ has $\lambda = 1$ if and only if
$w=\overline{2} 1^x 2^y 3^z$ for $x,y,z\geq 1.$
To make the braid closure one component, we
further require $x+z= \mbox{odd.}$
If $k, l \geq 1,$
$$w \leadsto w' \ \stackrel{\eqref{BB}}{\lra} \
\ol{2}\ (\ol{1}3) = 1 \ol{2}\ \ol{2} \Lra  1 \ol{2}$$  and $\lambda \geq 2.$

\noindent\textbf{Case 3-v$'$}\qua
When $k=l=0,$
$$w \leadsto w'= \ol{2}1.$$
Thus $w$ has $\lambda=1$ if and only if $w=\ol{2}\ \ol{2}1^x,$ which has $2$ or
$3$ components.
When $k=0$ and $l\geq 1$,
$$w \leadsto w'= \ol{2} (123)^l 1 \ \stackrel{\eqref{S}}{\lra} \
1 \ol{2}(123)^{l-1}1  \ \stackrel{\eqref{A}}{\lra} \ \cdots
 \ \stackrel{\eqref{A}}{\lra} \ 1 \ol{2}1 = 11\ol{2} \lra 1 \ol{2}$$
thus $w$ has $\lambda=1$ if and only if  $w=\ol{2}\ 1^x\ 2^y\ 3^z\ 1^w$ for
some $x,y,z,w\geq 1.$
When $k\geq 1, l=0$,
$$w \leadsto w' = \ol{2}\ (\ol{1}\ \ol{3}\ \ol{2})^k 1 =
1 \ol{2}\ (\ol{1}\ \ol{3}\ \ol{2})^k \ \stackrel{\eqref{W}}{\lra} \
\cdots  \ \stackrel{\eqref{W}}{\lra} \ 1 \ol{2}$$ and $\lambda \geq 2.$
When $k, l\geq 1$
$$w \leadsto w' \ \stackrel{\eqref{BB}}{\lra} \ \ol{2}(\ol{1}3)1 \Lra 1\ol{2}1
\lra 1\ol{2}$$ and $\lambda \geq 2.$

\noindent\textbf{Case 3-viii$'$}\qua
When $k=l=0,$
$$w \leadsto w'=\ol{3}\ \ol{2}1 = 1\ol{3}\ \ol{2}=\ol{2}1\ol{2}\Lra 1\ol{2}$$
Thus $w$ has $\lambda=1$ if and only if $w=\ol{3}\ \ol{2}1^x$ for
some $x\geq 1.$
When $k=0$ and $l\geq 1$,
$$w \leadsto w'=\ol{3}\ \ol{2}(123)^l 1 = 1\ol{3}\ \ol{2}(123)^l \Lra
1 \ol{2}(123)^l  \ \stackrel{\eqref{A}}{\lra} \
\cdots  \ \stackrel{\eqref{A}}{\lra} \ 1 \ol{2}$$ and $\lambda \geq 2.$
When $k\geq 1$ and $l=0$,
$$w \leadsto w'=\ol{3}\ \ol{2}(\ol{1}\ \ol{3}\ \ol{2})^k 1 \Lra
1 \ol{2}(\ol{1}\ \ol{3}\ \ol{2})^k  \ \stackrel{\eqref{W}}{\lra} \
\cdots  \ \stackrel{\eqref{W}}{\lra} \ 1 \ol{2}$$ and $\lambda \geq 2.$
When $k, l\geq 1$
$$w \leadsto w' \ \stackrel{\eqref{BB}}{\lra}\ \ol{3}\ \ol{2}(\ol{1}3)1 \Lra
1 \ol{2}\ \ol{1}3 = \ol{2} 33\ol{2} \Lra\ \lra  \ol{2}3 = 1 \ol{2}.$$
\newpage
and $\lambda \geq 2.$

\noindent\textbf{Case 3-ix$'$}\qua
When $k=l=0,$
$$w \leadsto w'=\ol{3}\ \ol{2}12 = 2 \ol{3}\ \ol{3}2 \Lra\ \lra 2 \ol{3}.$$
 Thus $w$ has
$\lambda=1$ if and only if $w= \ol{3}\ \ol{2}1^x 2^y$ for some $x,y \geq 1.$
When $k=0, l\geq 1$,
\begin{eqnarray*}
w &\leadsto &
w'=\ol{3}\ \ol{2}(123)^l 12=12\ol{3}\ \ol{2}(123)^l=\ol{2}11\ol{2}
(123)^l\lra\ \Lra \ol{2}3(123)^l  \\
&=& 1\ol{2}(123)^l\ \stackrel{\eqref{A}}{\lra} \
\cdots  \ \stackrel{\eqref{A}}{\lra} \ 1 \ol{2}
\end{eqnarray*}
and $\lambda \geq 2.$
When $k\geq 1, l=0$,
$$
w \leadsto w'=\ol{3}\ \ol{2}(\ol{1}\ \ol{3}\ \ol{2})^k12\lra\ \Lra 1\ol{2}
(\ol{1}\ \ol{3}\ \ol{2})^k\ \stackrel{\eqref{W}}{\lra} \
\cdots  \ \stackrel{\eqref{W}}{\lra} \ 1 \ol{2}$$ and $\lambda \geq 2.$
When $k, l\geq 1$
$$w \leadsto w'\ \stackrel{\eqref{BB}}{\lra}\ \ol{3}\ \ol{2}(\ol{1}3)12\lra\
\Lra 1\ol{2}(\ol{1}3)=11\ol{2}\ \ol{2}\lra\ \Lra 1\ol{2}$$ and
$\lambda \geq 2.$

\fullref{lambda-table} summarizes all the words with $\lambda=1.$

\begin{table}[ht!]
 \caption{}
 \label{lambda-table}
  \begin{center}
   \begin{tabular}{|l|l|}
\hline case & word with $\lambda =1.$ \\ \hline %
\strutt i & none. \\ \hline %
\strutt ii & $\overline{2}\ \overline{1}\
\overline{3}\ 1^{x}\ $(2 or 3 components.) \\ \hline %
v &
\begin{tabular}{l}
$\overline{3}\ 1^{x}\ 2^{y}\ 3^{z}\ 1^{w}\ =\left\{%
 \begin{array}{ll}
  \strutt C_{x+1,y,z} &\ \mbox{when } w=1, \\
  D_{x+1,y,z,w-1} & \ \mbox{when } w\geq 2.
 \end{array} \right.
$ \\
\strutt $\overline{3}\ \overline{3}\ 1^{x}\ $(2 or 3 components.)%
\end{tabular}
\\ \hline%
vi &
\begin{tabular}{l}
\strutt $\overline{3}\ \overline{3}\ 1^{x}\ 2^{y}=:B_{x,y}.$ \\ %
$\overline{3}\ 1^{x}\ 2^{y}\ 3^{z}\ 1^{w}\ 2^{v}\ =\left\{%
\begin{array}{ll}
C_{x+v+1,y,z} & \  \mbox{when } w=1, \\
\strutt D_{x+v+1,y,z,w-1} &\  \mbox{when } w\geq 2.%
\end{array} \right.
$%
\end{tabular}
\\ \hline %
\strutt vii & $\overline{1}\ \overline{3}\ 1^{x}\ 2^{y}\ 3^{z}\ =%
\overline{1}\ \overline{3}\ 1^{x+z}\ 2^{y}$ \\ \hline %
ix & $\overline{1}\ \overline{3}\ 1^{x}\ 2^{y}\ = \left\{
\begin{array}{lll}
\strutt \ol{3}\ \ol{3}1^{x+1} & \ \mbox{when } y=1 & (\mbox{$2$ or $3$ components)} \\
\strutt B_{x+1,y-1} & \ \mbox{when } y\geq 2 &
\end{array}\right. $
\\ \hline  %
\strutt i$'$ & none. \\ \hline %
\strutt iii$'$ & none. \\ \hline
\strutt iv$'$ & $\overline{2}\ 1^{x}\ 2^{y}\ 3^{z}\ =:C_{x,y,z}.$ \\
\hline v$'$ &
\begin{tabular}{l}
$\strutt \overline{2}\ 1^{x}\ 2^{y}\ 3^{z}\ 1^{w}\
=:D_{x,y,z,w}. $ \\
$\strutt \overline{2}\ \overline{2}\ 1^{x}\ $(2 or 3 components.)%
\end{tabular}\\ \hline %
\strutt viii$'$ & $\overline{3}\ \overline{2}\ 1^{x}\
=:A_{x}.$ \\ \hline
ix$'$ & $\overline{3}\ \overline{2}\ 1^{x}\ 2^{y}$ $=\left\{
\begin{array}{lll}
\strutt \overline{3}\ \overline{3}\ 2^{y+1} & \  \mbox{when } x=1, & (\mbox{$2$ or $3$
components}) \\
\strutt B_{x-1,y+1} &\ \mbox{when } x\geq 2. &
\end{array} \right.
$ \\ \hline
   \end{tabular}
   \end{center}
\end{table}
Words $A_x, \cdots, D_{x,y,z,w}$ are defined in \fullref{lambda-table}.
We can see that any
word with $\lambda=1$ and having one component has one of the forms;
$A_x, \cdots, D_{x,y,z,w}.$
\end{proof}

\begin{lemma}\label{alexander}
The leading terms of the Alexander polynomials of ${\cal K}_{n}$, $A_x$,
$B_{x, y}$, $C_{x, y, z}$ and $D_{x, y, z, w}$ are the following:
\begin{eqnarray*}
{\cal K}_n; && \pm(1 - 4t -6t^2 + 8t^3 - \cdots) \quad \mbox{if $n \geq 2$,} \\
A_x;        && \pm(1 - 3t + \cdots) \quad \mbox{if } x\geq 2, \\%
B_{x, y};   && \pm(1 - 3t + \cdots) \quad \mbox{if } x, y \geq 3, \\%
C_{x, y, z}; && \pm(1 - 5t + \cdots) \quad \mbox{if $x,y,z \geq 2$,}\\ %
C_{1, 2, z},  C_{1, y, 2}, C_{2, y, 1}, C_{x, 2, 1}; %
            && \pm(1 - 4t + 6t^2 -7t^3 + \cdots) %
            \quad \mbox{if $x, y, z \geq 4$,}\\ %
C_{1, y, z}, C_{x, y, 1}; %
            && \pm(1 - 4t + 7t^2 + \cdots) \quad \mbox{if $x, y, z \geq 3$}\\%
D_{x, y, z, w}, D_{x, y, z, 1}; && \pm(1 - 6t + \cdots) \quad \mbox{if } x,y,z,w \geq 2, \\
D_{x, y, 1, w}; && \pm(1 - 5t + \cdots) \quad \mbox{if } x,y,w \geq 2. %
\end{eqnarray*}
In particular, ${\cal K}_n\neq A_x, B_{x, y},C_{x, y, z}, D_{x, y, z, w}.$
\end{lemma}

\begin{figure}[hb!]
\begin{center}
\labellist\small\hair1.5pt
\pinlabel {(1)} at -15 425 
\pinlabel {(2)} at 358 348
\pinlabel {(3)} at -5 360   
\pinlabel {(4)} at 217 425  
\pinlabel {(5)} at 421 390 
\pinlabel {$x$} [t] at 44 304 
\pinlabel {$y$} [t] at 256 366 
\pinlabel {$z$} [t] at 464 300
\endlabellist
\includegraphics [height=40mm]{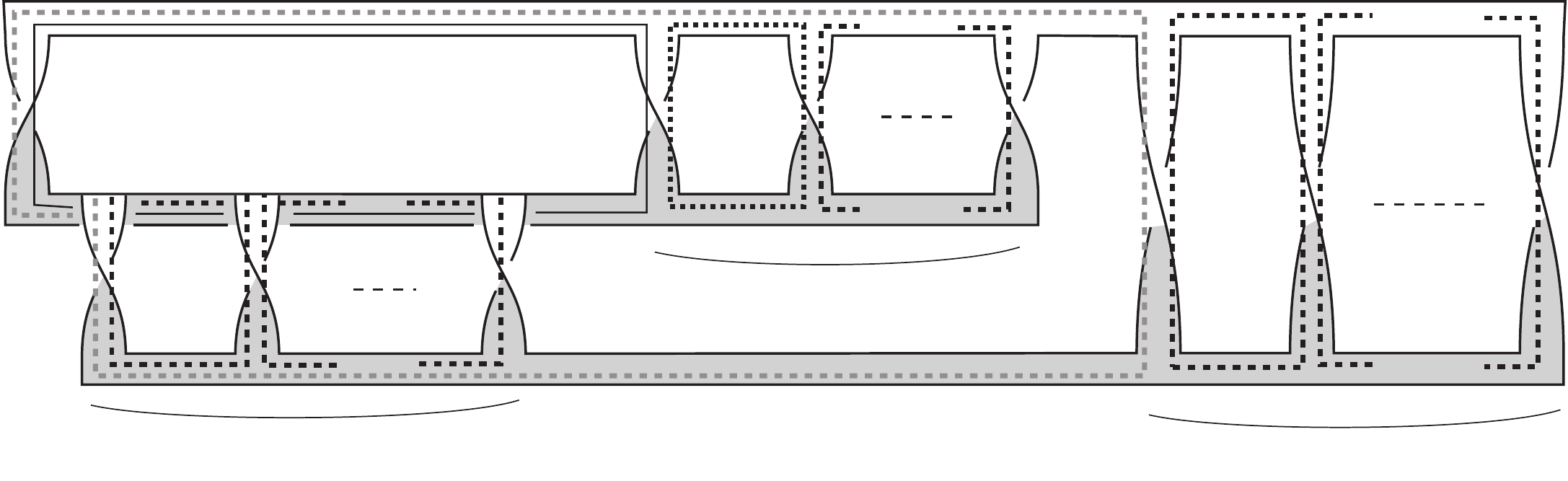}
\end{center}
\caption{The Bennequin surface $F$ of $C_{x, y, z}=\ol{2}
1^x 2^y 3^z$ and a basis for $H_1(F)$} \label{2-123}%
\end{figure}

\newpage
\begin{proof}[Proof of \fullref{alexander}] 
We prove that the
Alexander polynomial of $C_{x, y, z}$ for $x, y, z \geq 2$ is
$\pm(1 - 5t + \cdots)$.  Recall that Xu's Bennequin surface
is a minimal genus Seifert surface. Let $F$ be the
Bennequin surface of $C_{x, y, z}$  and choose a basis
$$\{ u^{(1)}, u^{(2)}, u^{(3)}_1, \cdots, u^{(3)}_{x-1}, u^{(4)}_1,
\cdots, u^{(4)}_{y-1}, u^{(5)}_1, \cdots, u^{(5)}_{z-1} \}$$ for
$H_1(F)$ as in \fullref{2-123},
where $u^{(k)}$ ($k=1,\cdots,5$) corresponds to loop $(k)$.

With respect to the basis, let $V_{x,y,z}$ denote the Seifert
matrix for $C_{x, y, z}$.
$$ V_{x,y,z} = \left[
  \begin{array}{c|c|cccc|cccc|cccc}
& \scriptstyle{1} &  &  &  &  & \scriptstyle{1} &  &  &  &  &  &  &  \\ \hline %
\scriptstyle{1} &  & \scriptstyle{-1} &  & &  &  &  &  &  & \scriptstyle{1} &  &  &  \\ \hline%
&  & \scriptstyle{-1} & \scriptstyle{1} &  &  &  &  &  &  &  &  &  &  \\
&  &  & \scriptstyle{-1} & \scriptstyle{\ddots}  &  &  &  &  &  &  &  &  &  \\
&  &  &  & \scriptstyle{\ddots}  & \scriptstyle{1} &  &  &  &  &  &  &  &  \\
&  &  &  &  & \scriptstyle{-1} &  &  &  &  &  &  &  &  \\ \hline%
&  &  &  &  &  & \scriptstyle{-1} & \scriptstyle{1} &  &  &  &  &  &  \\
&  &  &  &  &  &  & \scriptstyle{-1} & \scriptstyle{\ddots}  &  &  &  &  &  \\
&  &  &  &  &  &  &  & \scriptstyle{\ddots}  & \scriptstyle{1} &  &  &  &  \\
&  &  &  &  &  &  &  &  & \scriptstyle{-1} &  &  &  &  \\ \hline%
&  &  &  &  &  &  &  &  &  & \scriptstyle{-1} & \scriptstyle{1} &  &  \\
&  &  &  &  &  &  &  &  &  &  & \scriptstyle{-1} & \scriptstyle{\ddots}  &  \\
&  &  &  &  &  &  &  &  &  &  &  & \scriptstyle{\ddots}  & \scriptstyle{1} \\
&  &  &  &  &  &  &  &  &  &  &  &  & \scriptstyle{-1}%
\end{array} \right] $$
The empty spaces contain only $0$'s. The $3$rd (resp.\ $4$th,
$5$th) diagonal block has size $(x-1)\times (x-1)$ (resp.\
$(y-1)\times (y-1)$, $(z-1)\times (z-1)$). The Alexander polynomial
satisfies:

$\Delta_{x,y,z}(t) = \det(V_{x,y,z}^T - tV_{x,y,z})$
$$
= \det \left[
\begin{array}{c|c|cccc|cccc|cccc}
& \scriptstyle{1-t} &  &  &  &  & -t &  &  &  &  &  &  &  \\ \hline%
 \scriptstyle{1-t} &  & t &  &  &  &  &  &  &  & -t &  &  &  \\ \hline%
& -1 & \scriptstyle{-1+t} & -t &  &  &  &  &  &  &  &  &  &  \\
&  & 1 & \ddots  & \ddots  &  &  &  &  &  &  &  &  &  \\
&  &  & \ddots  & \ddots  & -t &  &  &  &  &  &  &  &  \\
&  &  &  & 1 & \scriptstyle{-1+t} &  &  &  &  &  &  &  &  \\ \hline%
1 &  &  &  &  &  & \scriptstyle{-1+t} & -t &  &  &  &  &  &  \\
&  &  &  &  &  & 1 & \ddots  & \ddots  &  &  &  &  &  \\
&  &  &  &  &  &  & \ddots  & \ddots  & -t &  &  &  &  \\
&  &  &  &  &  &  &  & 1 & \scriptstyle{-1+t} &  &  &  &  \\ \hline%
& 1 &  &  &  &  &  &  &  &  & \scriptstyle{-1+t} & -t &  &  \\
&  &  &  &  &  &  &  &  &  & 1 & \ddots  & \ddots  &  \\
&  &  &  &  &  &  &  &  &  &  & \ddots  & \ddots  & -t \\
&  &  &  &  &  &  &  &  &  &  &  & 1 & \scriptstyle{-1+t}%
\end{array} \right]. $$
Expanding it by the $(x+1)$th column, we have;

$ \Delta_{x,y,z}(t) = (-1+t)\Delta_{x-1,y,z}(t)$
$$ -(-t)\det
\left[
\begin{array}{c|c|cccc|ccc|ccc}
& \scriptstyle{1-t} &  &  &  &  & -t &  & &  &  &  \\ \hline%
\scriptstyle{1-t} &  & t &  &  &  &  &  &  & -t &  &  \\ \hline%
& -1 & \scriptstyle{-1+t} & -t &  &  &  &  &  &  &  &  \\
&  & 1 & \ddots  & -t &  &  &  &  &  &  &  \\
&  &  & 1 & \scriptstyle{-1+t} & -t &  &  &  &  &  &  \\
&  &  &  &  & 1 &  &  &  &  &  &  \\ \hline%
1 &  &  &  &  &  & \scriptstyle{-1+t} & -t &  &  &  &  \\
&  &  &  &  &  & 1 & \ddots  & -t &  &  &  \\
&  &  &  &  &  &  & 1 & \scriptstyle{-1+t} &  &  &  \\ \hline%
& 1 &  &  &  &  &  &  &  & \scriptstyle{-1+t} & -t &  \\
&  &  &  &  &  &  &  &  & 1 & \ddots  & -t \\
&  &  &  &  &  &  &  &  &  & 1 & \scriptstyle{-1+t}%
\end{array} \right] $$
$= (-1+t)\Delta_{x-1,y,z}(t) + t \Delta_{x-2,y,z}(t).$

\eject If $\Delta_{i,y,z}(t)=(-1)^i (\alpha_0 + \alpha_1 t + \alpha_2 t^2
+ \cdots )$ for $i=x-1$ and $x-2$, then
\begin{eqnarray*}
\Delta_{x,y,z}(t) &=& (-1+t) (-1)^{x-1} (\alpha_0 + \alpha_1 t +
\alpha_2 t^2 + \cdots )\\
&&\hspace{2in}+ t (-1)^{x-2} (\alpha_0 + \alpha_1 t + \alpha_2 t^2 + \cdots) \\ %
&=& (-1)^x (\alpha_0 + \alpha_1 t + \alpha_2 t^2 + \cdots ).
\end{eqnarray*}
In fact, $\Delta_{x,y,z}(t) = (-1)^{x+y+z} (1-5t+ \cdots )$ for
all $x,y,z \in \{2,3\}$.  By induction, $\Delta_{x,y,z}(t) =
(-1)^{x+y+z} (1-5t+ \cdots )$ for all $x, y, z \geq 2.$

Other cases follow by similar arguments.
\end{proof}

Finally we are ready to prove the theorem.

\begin{proof}[Proof of \fullref{BM-thm}] 
By Lemmas
\ref{lambda=1}, \ref{ABCD}, \ref{alexander}, our knot ${\cal
K}_{2m}$ where ($m\geq 1$) cannot be a $3$--braid. Thus by \fullref{D+}, \fullref{BM-thm} follows.
\end{proof}

\section{Uniqueness of the
algebraic crossing number at minimal braid index}\label{chap3}

\subsection{Sharpness of the MFW--inequality and conjectures}

It has been conjectured (see \cite[page 357]{Jones-1} for
example) that the exponent sum in a minimal braid representation is
a knot invariant.
\begin{conjecture}[Main Conjecture]\label{Jones-conj}
Let ${\cal K}$ be a knot type of braid
index $b_{\cal K}$.  If $K^1$ and $K^2$ are braid representatives of
${\cal K}$ with $b_{K^1}=b_{K^2}=b_{\cal K}$ then their algebraic crossing numbers
have $c_{K^1}=c_{K^2}$.
\end{conjecture}

We deform it into: 

\begin{conjecture}[Stronger Conjecture]\label{greedy}
Let ${\cal B}_{\cal K}$ be the set of braid representatives of ${\cal K}.$
Let $\Phi\co {\cal B}_{\cal K} \to \mathbb{N} \times \mathbb{Z}$ be a map such
that  $\Phi(K):=\left( b_K, c_K \right)$ for $K\in {\cal B}_{\cal K}.$
Then there exists a unique $c_{\cal K} \in \mathbb{Z}$ with
\begin{equation}\label{quadrant equation}
\Phi({\cal B}_{\cal K})=\left\{ (b_{\cal K}+x+y, c_{\cal K}+x-y)\ |\
x,y\in \mathbb{N} \right\},
\end{equation}
a subset of the infinite quadrant region shaded in \fullref{b-region}.
\end{conjecture}
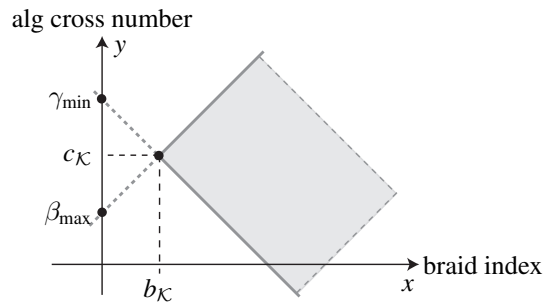
\begin{figure}[ht!]
\begin{center}
\labellist\small\hair2pt
\pinlabel {$\gamma_{\min}$} [r] at 49 147
\pinlabel {$\beta_{\max}$} [r] at  49 68
\pinlabel {$c_{\cal K}$} [r] at 49 108
\pinlabel {$b_{\cal K}$} [t] at 92 25
\pinlabel {alg cross number} [b] at 52 192
\pinlabel {braid index} [l] at 275 31
\pinlabel {$x$} [t] at  268 25
\pinlabel {$y$} [l] at 57 184
\endlabellist
\includegraphics[height=40mm]{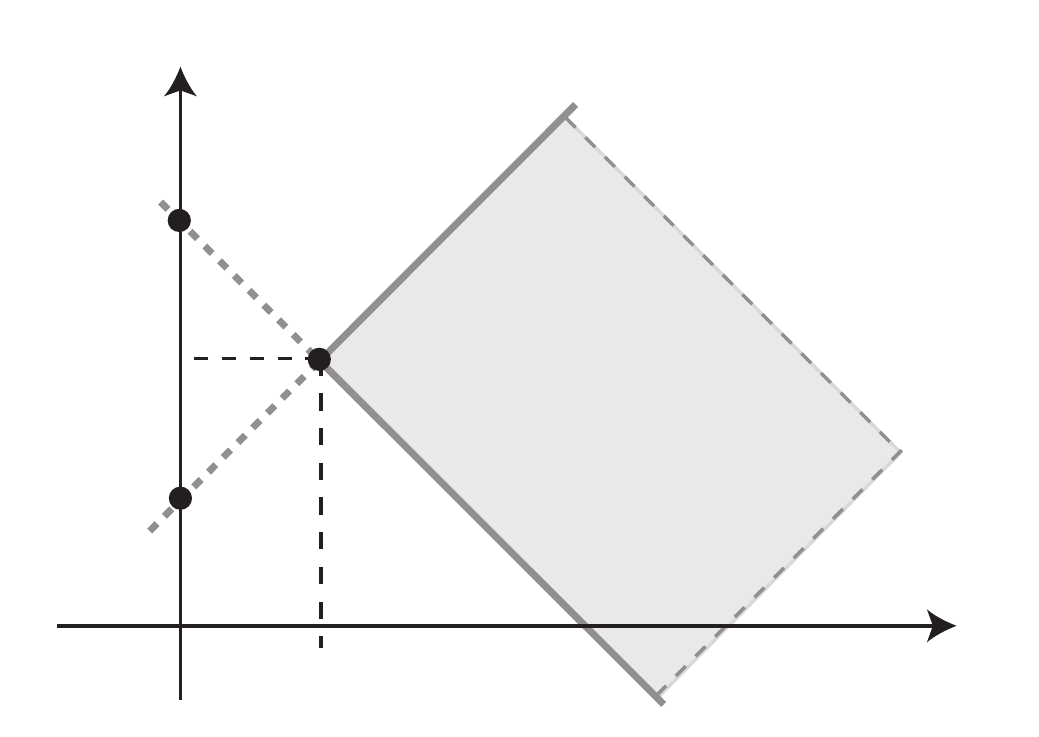}
\caption{The region of braid representatives of ${\cal K}$}\label{b-region}
\end{center}
\end{figure}

The inclusion ``$\supset$'' is trivial by the following argument: Let $K_\star
\in B_{\cal K}$ be a minimal braid representative with
$\Phi(K_\star)=(b_{\cal K}, c_{\cal K}).$ Suppose $K\in
B_{\cal K}$ is obtained from $K_\star$ after applying $(+)$--stabilization $x$--times
and then $(-)$--stabilization $y$--times.  Then $(b_{\cal K}+x+y,
c_{\cal K}+x-y)=\Phi(K) \in \Phi({\cal B}_{\cal K}).$

The MFW--inequality \eqref{lower bound-1} says that
$c_K \geq -b_K+(d_++1),\ c_K\leq b_K+(d_--1)$
for any $K\in {\cal B}_{\cal K}.$
Thus
\begin{eqnarray}\label{quad-eq}
\Phi({\cal B}_{\cal K})\subset \left\{(x,y)\ |\ b_{\cal K}\leq x,\
-x+(d_++1)\leq y \leq x+(d_--1)\right\}.
\end{eqnarray}

Before we provide examples of the conjectures we present:

\begin{theorem}\label{imply}
Sharpness of the MFW--inequality implies the truth of Conjectures
\ref{Jones-conj} and \ref{greedy}. In particular;
$$b_{\cal K}=\frac{d_+-d_-}{2}-1, \quad c_{\cal K}=\frac{d_++d_-}{2}.$$
\end{theorem}
We remark that 
the statement in the theorem with regard to \fullref{Jones-conj}
has been well known to many people.

\begin{proof}[Proof of \fullref{imply}]
Let $K_\star
\in {\cal B}_{\cal K}$ be a minimal braid representative. Since the MFW
inequality \eqref{lower bound-1} is sharp on ${\cal K},$ we have
$c_{K_\star}-b_{\cal K}+1=d_-,$ and $d_+=b_{\cal K}+c_{K_\star}-1,$
ie, $c_{K_\star}=(d_++d_-)/2$ which
is independent of the choice of $K_\star.$ Thus we denote $c_{K_\star}=:c_{\cal K}.$
In this case, the right side of \eqref{quad-eq} coincides with the right side
of \eqref{quadrant equation} and we have the other inclusion
 ``$\subset$'' of \eqref{quadrant equation}.
\end{proof}

\begin{example}\label{ex-of-conj}{\rm
Both of the conjectures are true for unlinks,
torus links, closed positive braids with
a full twist (for example, the Lorenz links) \cite{FW},
$2$--bridge links and alternating fibered links \cite{Murasugi}, where the
MFW--inequality is sharp and one can apply \fullref{imply}.

Also \fullref{Jones-conj} applies to
links with braid index $\leq 3$ \cite{BM3}.
However, this case has been settled by a completely different
way, the classification of $3$--braids. Namely, any link of braid index $3$
admits a unique conjugacy class of $3$--braid representatives or has at most two
conjugacy classes of $3$--braid representatives related to each other by a
flype move, which does not change the algebraic crossing number of the link.
}\end{example}

Every transversal knot $TK$ in $S^3$ with the standard contact
structure is transversally isotopic to a transversal closed braid
$K$ \cite{Ben}. The {\it Bennequin number} $\beta$ is an invariant
of transversal knots. By the identification of $TK$ and $K$, we
have $\beta (K) = c_K - b_K.$ If \fullref{greedy} is true for
${\cal K},$ then the
{\it maximal} Bennequin number $\beta_{\rm max}({\cal K})$ for the knot
type ${\cal K}$ is realized on ${\cal B}_{\cal K} \ni K$'s  with
plotted vertices  $\Phi(K)$ on the upper half boundary of the quadrant
region of \fullref{b-region}. Let
$\gamma (K) := c_K + b_K.$
For any $K$ and its mirror image $\overline{K},$ 
we have $\gamma(K)=-\beta(\overline{K}).$
Thus $\gamma_{\rm min}({\cal K}) = - \beta_{\rm max}
(\overline{\cal K}).$
See \fullref{b-region}. Thus, investigation of
$\beta_{\rm max}({\cal K})$ is related to \fullref{greedy}.

\subsection{Cabling and the conjectures}
In this subsection, we study behavior of the deficit of the MFW
inequality under cabling and prove \fullref{deficit-cable}. As
a consequence, we observe that the \fullref{Jones-conj} is
true for many of the knots and links that appeared in \fullref{chap2},
where the MFW--inequality is not sharp, ie, we cannot apply
\fullref{imply}.
\vspace{-2pt}

We also prove, in Theorems \ref{cable-thm} and \ref{cable-thm-link},
that the truth of \fullref{Jones-conj} is
``inherited'' through cabling operations.
\vspace{-2pt}

\begin{figure}[ht!]
\begin{center}
\labellist\small
\pinlabel {$k$--times} [t] at 295 396
\endlabellist
\includegraphics [height=40mm]{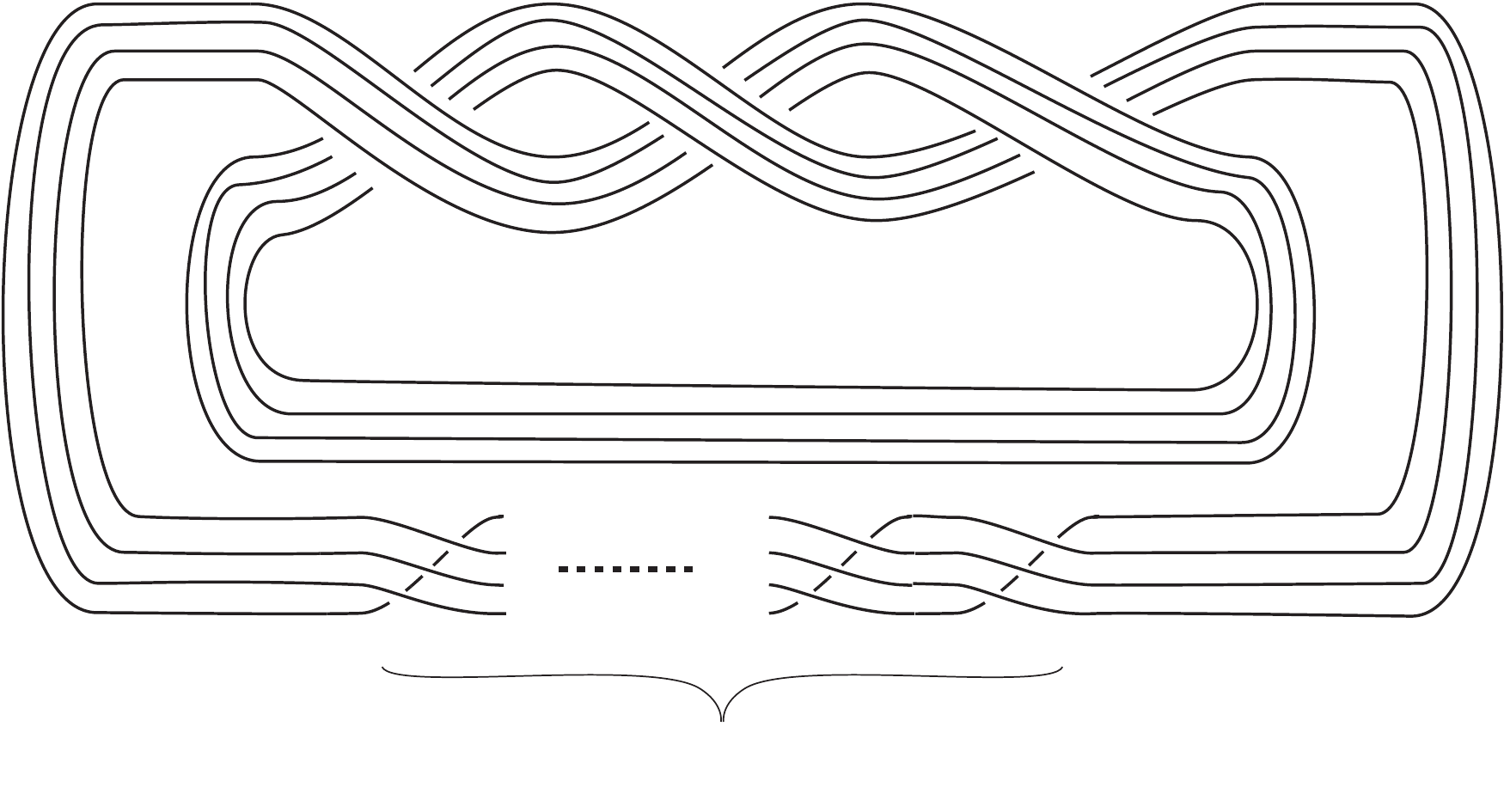}
\end{center}
\caption{$(4,q)$--cable ($q=4 \cdot 3 + k$) of the right hand trefoil}
\label{(4,q)-trefoil}%
\end{figure}

Let us fix some notation.  Let ${\cal K}$ be a knot type.
Denote the $(p,q)$--cable of ${\cal K}$ by ${\cal K}_{p,q}.$ Let
$K$ be a braid representative with $b_K = b_{\cal K}$ and with
algebraic crossing number $c_K$. Put $$k:= q - p\cdot c_K$$ and let
$K_{p,q}$ denote the $p$--parallel copies of $K$ with a $k/p$--twist
(see \fullref{(4,q)-trefoil}).
We can assume that $K_{p,q}$ is on the boundary of a
tubular neighborhood $N$ of $K$ (thus $K$ is the core of solid
torus $N$). Then
$$q={\rm lk}(K_{p,q}, K)$$ where `lk' is the linking number.
\fullref{(4,q)-trefoil} shows that the algebraic crossing number of $K_{p,q}$ is
\begin{equation}\label{c_pq}
c_{K_{p,q}}=p^2 c_K + k(p-1)=q(p-1)+p\cdot c_K.
\end{equation}
Thanks to \cite{W}, we know that the braid index of ${\cal K}_{p,q}$ satisfies
\begin{equation}\label{cable-index}
b_{{\cal K}_{p,q}}=p\cdot b_{\cal K}.
\end{equation}

Although one can see a similar result in Theorem 7 of \cite{S}, we state the
following for completeness of our discussion:
\begin{theorem}\label{deficit-cable}
Suppose $K^1$ and $K^2$ are braid representatives of
${\cal K}$ with $b_{K^1}=b_{K^2}=b_{\cal K}$
and with distinct algebraic crossing numbers
$c_{K^1} < c_{K^2},$ (ie, \fullref{Jones-conj} does not apply to
${\cal K}).$ Then the deficit $D_{{\cal K}_{p,q}}$
of the MFW--inequality for $(p,q)$--cable ${\cal K}_{p,q}$ is;
\begin{equation}\label{eq-deficit-cable}
D_{{\cal K}_{p,q}} \geq \frac{p}{2} ( c_{K^2}-c_{K^1} ) \geq p.
\end{equation}
\end{theorem}
\begin{proof}[Proof of \fullref{deficit-cable}]
Thanks to \eqref{cable-index}, and by the construction of
$({K^1})_{p,q}$ and $({K^2})_{p,q}$, they are both
minimal braid representatives of ${\cal K}_{p,q}$ ie,
$b_{{K^1}_{p,q}}=b_{{K^2}_{p,q}}=b_{{\cal K}_{p,q}}=p\cdot b_{\cal K}.$

Let $k_1, k_2$ be integers satisfying
$q=p c_{K^1} + k_1 = p c_{K^2} + k_2.$
By \eqref{c_pq} we have
$
c_{{K^1}_{p,q}} = p^2 c_{K^1} + k_1 (p-1)$ and
$c_{{K^2}_{p,q}} = p^2 c_{K^2} + k_2 (p-1).$
Therefore,
$c_{{K^2}_{p,q}}-c_{{K^1}_{p,q}} = p(c_{K^2}-c_{K^1}).$
By \eqref{lower bound-1} we have
$c_{{K^2}_{p,q}} - b_{{\cal K}_{p,q}} +1 \leq d_{-}\leq d_{+}\leq
 c_{{K^1}_{p,q}} + b_{{\cal K}_{p,q}} -1,$
and by \fullref{deficit-def},
\begin{equation}\label{D(cable)}
D_{{\cal K}_{p,q}} \geq \frac{1}{2}  (c_{{K^2}_{p,q}}-c_{{K^1}_{p,q}}) =
\frac{p}{2}  (c_{K^2}-c_{K^1}).
\end{equation}
This is the first inequality of \eqref{eq-deficit-cable}.

Notice that $K^1$ and $K^2$ are related each other by a sequence of Markov
moves \cite{B}. Let $K^1 = B_1 \rightarrow B_2 \rightarrow \cdots
\rightarrow B_n=K^2$ be a Markov tower.  Each arrow corresponds to
either braid isotopy, stabilization or destabilization moves. Let
$(x_i, y_i)$ be the braid index and the algebraic crossing number of
$B_i.$ Then  $(x_{i+1}, y_{i+1})-(x_i, y_i) =  (0,0), (\pm 1, \pm
1)$ or $(\mp 1, \pm 1)$ depending on the move corresponding to the
arrow between $B_{i+1}$ and $B_i$. Since $x_1 = x_n = b_{\cal K}$
the difference $c_{K^1}-c_{K^2}=y_1-y_n\neq 0$ must be an even integer.
Therefore, we get the second inequality of
\eqref{eq-deficit-cable}.
\end{proof}

\begin{corollary}\label{conj-5knots}
\fullref{Jones-conj} is true for all $9_{42}, 9_{49},
10_{132}, 10_{150}, 10_{156}.$
\end{corollary}

In \cite{K}, it is proved that \fullref{greedy} also holds for the five knots.

\begin{proof}
Knotscape computes that the deficit of $2$--cable
${\cal K}_{2,2c_K+1}$ is $1$ for each knot.
\begin{center}
\begin{tabular}{|l|c|c|c|c|c|}
\hline
\struttt ${\cal K}$ & $b_{\cal K}$ & $D_{\cal K}$ & $K$ & $c_K$ & $D_{{\cal K}_{2,2c_K+1}}$
\\ \hline
\strutt  $9_{42}$   & $4$ & $1$ & $aaacBAAcB$   & $1$ & $1$
\\ \hline
\strutt  $9_{49}$   & $4$ & $1$ & $aabbcbAbbcB$ & $7$ & $1$
\\ \hline
\strutt  $10_{132}$ & $4$ & $2$ & $AbcaaaBBBcb$ & $3$ & $1$
\\ \hline
\strutt  $10_{150}$ & $4$ & $1$ & $aabbcbABccB$ & $5$ & $1$
\\ \hline
\strutt  $10_{156}$ & $4$ & $1$ & $aaacBAAcbAb$ & $3$ & $1$
\\ \hline
\end{tabular}\end{center}
Comparing with \eqref{eq-deficit-cable}, each ${\cal K}$ must have unique
algebraic crossing number.
\end{proof}

Thanks to Knotscape, the $(2, 2c_{{\cal K}_n}+1)$--cable of
${\cal K}_n = BM_{-1, -2, n, 2}$ has deficit $=1$ if $|n|$ is small.
ie, \fullref{Jones-conj} is true for ${\cal K}_n$ if $|n|$ is small.
\fullref{conj-5knots} implies:

\begin{corollary}
\fullref{Jones-conj} is true for the prime links ${\cal A}^n(9_{42})$
(see \fullref{942}$).$
\end{corollary}

\begin{proof}
We know that $9_{42}$ has unique algebraic crossing number $=1$ by
\fullref{conj-5knots}. Since each link component of 
${\cal A}^n(9_{42})$ is
$9_{42},$ we get this corollary.
\end{proof}

With regard to the deficit of cabled links, we conjecture that:

\begin{conjecture}\label{limit}
For any $q$, the limit
$\displaystyle\lim_{p \rightarrow \infty}D_{{\cal K}_{p,q}}$ of deficits
exists.
\end{conjecture}

\begin{remark}{\rm
If \fullref{limit} is true, then \eqref{eq-deficit-cable} of
\fullref{deficit-cable} implies
the truth of \fullref{Jones-conj}.}
\end{remark}

We present another property of cabling:
\begin{theorem}\label{cable-thm}
Let ${\cal K}$ be a non-trivial knot type.
If \fullref{Jones-conj} is true for ${\cal K}$ then it is also
true for ${\cal K}_{p,q}$ when $p\geq 2.$

In particular, if $c_{\cal K}$ and $c_{{\cal K}_{p,q}}$ denote the unique
algebraic
crossing numbers of ${\cal K}$ and ${\cal K}_{p,q}$ respectively in their
minimal braid representatives then we have
$$c_{{\cal K}_{p,q}} =(p-1)q+p\cdot c_{\cal K}.$$
\end{theorem}

\begin{remark} {\rm Suppose ${\cal K}$ is the right hand trefoil. The
MFW--inequality is sharp on ${\cal K}.$ Since its cable ${\cal
K}_{2,7}$ has deficit $D_{{\cal K}_{2,7}}=1$ (see \cite{MS}), we cannot apply
\fullref{imply}. However \fullref{cable-thm}
guarantees the truth of the conjecture for ${\cal K}_{2,7}.$}
\end{remark}

The following proof is inspired by the work of Williams \cite{W}, whose main
result can be seen in formula \eqref{cable-index}. Note that his result holds
not only for cable knots but also for generalized cable links.
For the sake of completeness we repeat part of his discussion.

\begin{proof}[Proof of \fullref{cable-thm}]
Assume
\fullref{Jones-conj} is true for ${\cal K}$ and denote the
unique algebraic crossing number at minimal braid index by
$c_{\cal K}.$

Let $K$ be a braid representative of ${\cal
K}.$ Suppose $K'$ is a braid representative of ${\cal K}_{p,q}$ on
the boundary of a small tubular (solid torus) neighborhood $N$ of
$K.$ We may regard the $z$--axis as the braid axis. Let $\phi\co 
{\mathbb R}^3 \rightarrow {\mathbb R}^3$ be a diffeomorphism of
compact support so that 
\begin{itemize}
\item[$(\star)$] $\phi(K') \subset \partial\phi(N)$ has
exactly $p\cdot b_{\cal K}$ maxima and $p\cdot b_{\cal K}$ minima (both
non-degenerate critical points) and no other critical points and 
\item[$(\star\star)$]
the ``height'' function $h\co  \partial\phi(N)\simeq T^2 \rightarrow
\mathbb{R}$ is a Morse function. 
\end{itemize}
In particular, $\phi(K')$ has a
braid position with braid index $p\cdot b_{\cal K}.$

By $(\star\star)$, a generic intersection of the horizontal plane with
$\partial\phi(N)\simeq T^2$ consists of disjoint simple closed
curves.
Furthermore, these simple closed curves are either meridians of $T^2$
or trivial in $T^2$ since ${\cal K}$ is knotted (Remark 1 of \cite{W}).

Remark 2 in \cite{W} says that there is a plane $\pi$ (parallel to
the $(xz)$--plane) intersecting transversely with $T^2$ in a
meridian.

Let $J$ be an innermost one among such meridians. Then $J$ bounds
a disk $d\subset \pi \cap \phi(N)$ which separates $\phi(K')$ into arcs
$\{C_i\}$. Close each $C_i$ with aid of some arc $D_i \subset d$
and set $\hat{K_i}:=C_i \cup D_i.$ See \fullref{cable1}.

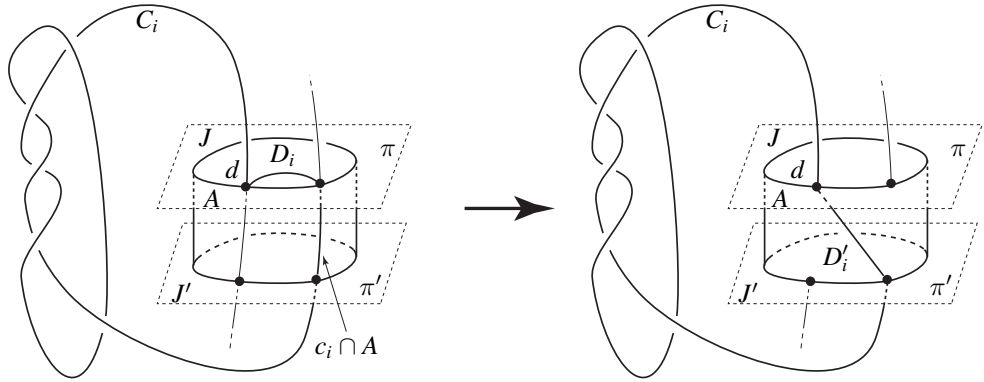
\begin{figure}[ht!]
\begin{center}
\labellist\small\hair1.5pt
\pinlabel {$C_i$} [t] at 105 551
\pinlabel {$C_i$} [t] at 445 551
\pinlabel {$D_i$} at 186 462
\pinlabel {$d$} at 156 454
\pinlabel {$d$} at 493 454
\pinlabel {$D_i'$} at 516 401
\pinlabel {$A$} [b] at 144 429
\pinlabel {$A$} [b] at 483 429
\pinlabel {$\pi$} [r] at 256 465
\pinlabel {$\pi$} [r] at 596 465
\pinlabel {$J$} [t] at 140 480
\pinlabel {$J$} [t] at 479 480
\pinlabel {$J'$} [b] at 126 372
\pinlabel {$J'$} [b] at 464 372
\pinlabel {$\pi'$} [r] at 248 385
\pinlabel {$\pi'$} [r] at 588 385
\pinlabel {$c_i \cap A$} [t] at 224 356
\endlabellist
\includegraphics [height=50mm]{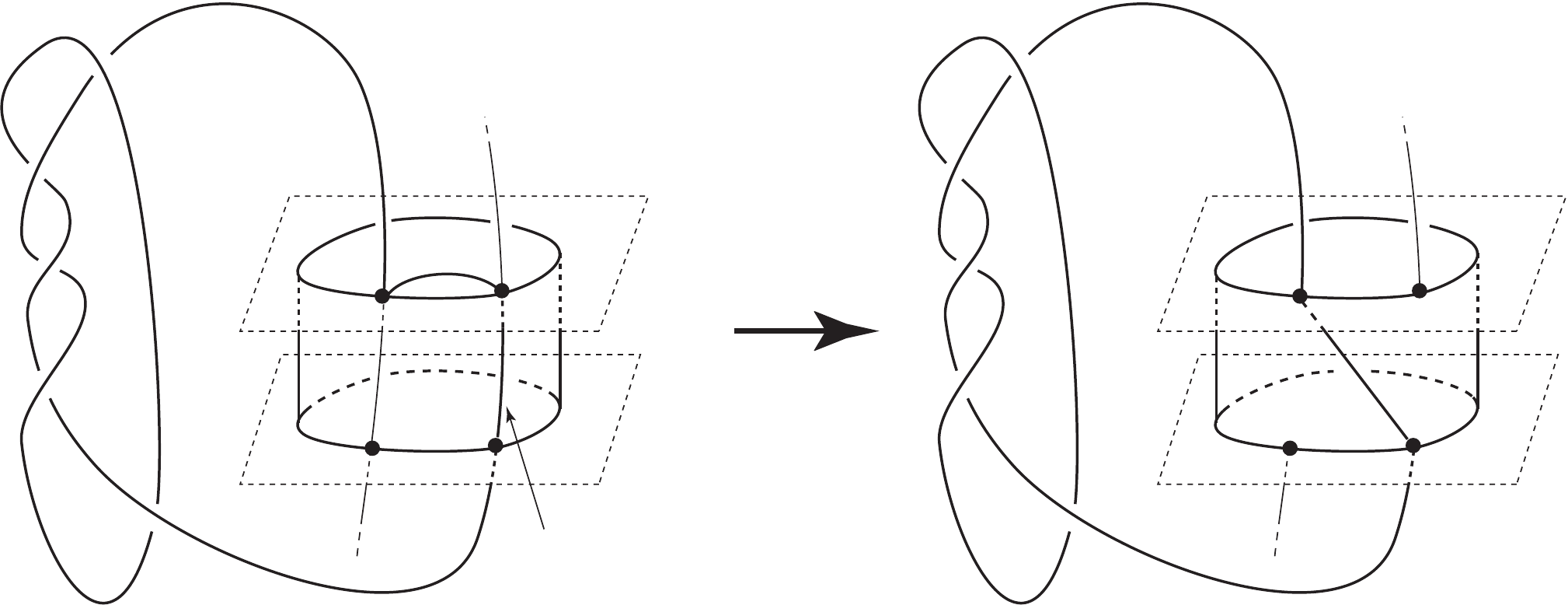}
\end{center}
\caption{Construction of $K_i$ from $\hat{K_i}$}\label{cable1}
\end{figure}

Thanks to Remark 3 in \cite{W}, $p$ of $\hat{K_i}$'s are
non-trivial (ie, do not bound any disk in $\phi(N))$ since the
linking number of  $J$ and $\phi(K')$ pushed a little bit into the
interior of $\phi(N)$ is $p.$

Discard trivial  $\hat{K_i}$'s.

Our $\hat{K_i}$'s are not in a braid position. As in \cite{W},
we make them have a braid position:
Choose another plane $\pi'$ just below $\pi$ and call the annulus
between the two planes $A$ (see \fullref{cable1}). We may
assume that the other boundary curve $J'\subset \partial A$ is
parallel to $J.$

As in the passage of \fullref{cable1} replace the arc $D_i \cup
(C_i \cap A)$ (the left sketch) with $D_i' \subset A$ (the right sketch)
and construct $p$--parallels;
\begin{equation}\label{K_i}
K_i:= (C_i - (C_i \cap A)) \cup D_i' \subset \partial\phi(N)\ \
\mbox{ for } i=1,\cdots,p,
\end{equation}
which is in a braid position. Also the $K_i$'s are disjoint from each other
and each is isotopic to the core of the solid torus
$\phi(K)\simeq {\cal K},$ thus $b_{\cal K} \leq
b_{K_i}.$ Then we have
\begin{eqnarray*}
p\cdot b_{\cal K} &\leq & \sum_{i=1}^p \{b_{K_i}= \mbox{number of max of $K_i$}\}\\
            &\leq & \{ \mbox{number of max of $\phi(K')$} \} = p\cdot b_{\cal K}
\end{eqnarray*}
where the last equality holds by $(\star\star)$ above. This implies that
\begin{itemize}
\item[$(\dag)$] there are no trivial $\hat{K_i}$'s (we didn't have to discard
anything), 
\item[$(\dag\dag)$] each knot has $b_{\cal K} = b_{K_i}.$ 
\end{itemize}
Let $n,
0\leq m <p$ be integers such that
$$q= p(c_{\cal K}+n)+m.$$
By $(\dag)$, the $p$--component link $L:= K_1 \cup \cdots \cup K_p$
is obtained from ${\bf K'}:= \phi(K')$ by using the meridian disk
$d$ to create a cutout and adding an $m/p$--twist along the annulus
$A,$ then gluing the end-points. See \fullref{twists}.
In other words, $L$ is the $(p, p(c_{\cal K}+n))$--cable
of $K$. From $(\dag\dag)$ we have
$c_{K_i}=c_K=c_{\cal K}.$ Therefore, $L$ has the algebraic crossing number
\begin{figure}[htpb!]
\begin{center}\includegraphics [height=30mm]{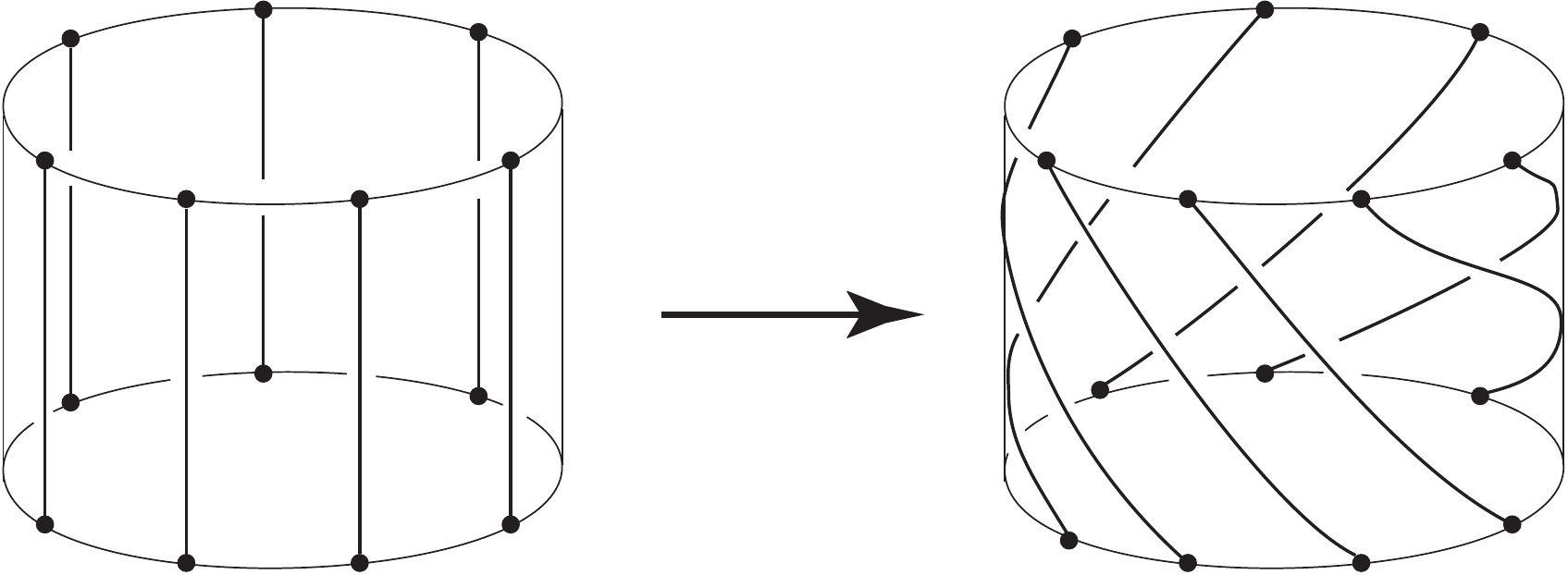}
\end{center}
\caption{From $A \cap {\bf K'}$ to $A \cap L$, where $p=7, m=2$}\label{twists}
\end{figure}
\newpage
\begin{eqnarray*}
c_L &=& \sum_{i=1}^p c_{K_i} + \sum_{i\neq j} {\rm lk}(K_i, K_j)
     = p\cdot c_{\cal K} + p(p-1)(c_{\cal K}+n) \\
    &=& p^2 c_{\cal K} + p(p-1)n
\end{eqnarray*}
and $c_{\bf K'}-c_L = m (p-1).$ Thus,
\begin{eqnarray*}
c_{\bf K'} &=&  p^2 c_{\cal K} + p(p-1)n +  m (p-1) \\
           &=&  p^2 c_{\cal K} + (p-1)(pn+m) \\
           &=&  p^2 c_{\cal K} + (q-p\cdot c_{\cal K}) \\
           &=&  (p-1)q+p\cdot c_{\cal K},
\end{eqnarray*}
which is independent of the choice of ${\bf K'}\in {\cal B}_{{\cal K}_{p,q}}.$
Compare with \eqref{c_pq}. This
concludes the uniqueness of the algebraic crossing number of
${\cal K}_{p,q}$ at minimal braid index.
\end{proof}

A similar result to \fullref{cable-thm} holds for links.

\begin{theorem}\label{cable-thm-link}
Let ${\cal L}={\cal K}^{(1)}\cup\cdots\cup{\cal K}^{(l)}$ be an $l$--component
link of braid index $=b_{\cal L}.$ Assume that each ${\cal K}^{(j)}$ is
a non-trivial knot.
Let ${\cal L}':={\cal K}^{(1)}_{p,q_1}\cup
\cdots\cup{\cal K}^{(l)}_{p,q_l}$ be the $p$--cable of ${\cal L}$
such that $q_j={\rm lk}({\cal K}^{(j)}, {\cal K}^{(j)}_{p,q_j})$ for $j=1,
\cdots,l.$

If ${\cal L}$ and every component ${\cal K}^{(j)}$ have unique algebraic
crossing numbers $c_{\cal L}, c_{{\cal K}^{(j)}}$
in minimal braid representations, then so does ${\cal L}'.$

Furthermore, let $k_j$ satisfy $q_j=p\cdot c_{{\cal K}^{(j)}}+k_j$ then
\begin{equation}\label{cable-thm-link-eq}
  c_{{\cal L}'}=p^2 c_{\cal L} + (p-1)(k_1+\cdots+k_l).
\end{equation}
\end{theorem}

\begin{remark}\label{flype-remark} {\rm
The assumption for ${\cal K}^{(j)}$ in the second paragraph of \fullref{cable-thm-link} is essential by the following reason:
If $b_{\cal L}=\sum_{j=1}^lb_{{\cal K}^{(j)}},$ then the existence
of unique algebraic crossing number of ${\cal L}$ in minimal braid
representation implies that each  ${\cal K}^{(j)}$ also has unique
algebraic crossing number.  However, this is not true in general. For instance,
assume that ${\cal L}={\cal K}^{(1)}\cup {\cal K}^{(2)}$ is a $2$--component
link and has two braid representatives related to each other by a flype move
as in Figure\ref{flype}.
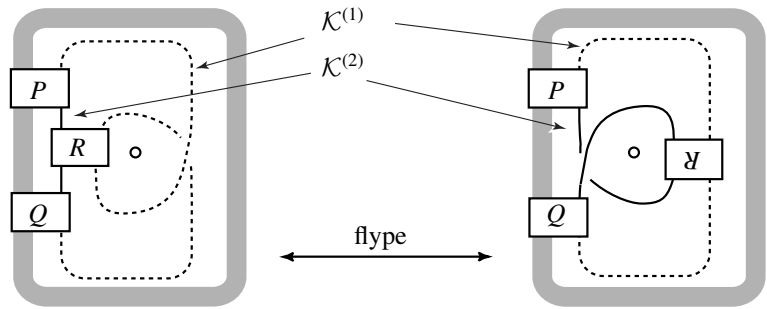
\begin{figure}[ht!]
\begin{center}
\labellist\small
\pinlabel {$P$} at 60 661
\pinlabel {$P$} at 296 661
\pinlabel {$Q$} at 60 604
\pinlabel {$Q$} at 296 604
\pinlabel {$R$} <0pt,1pt> at 78 634
\pinlabel \rotatebox{180}{$R$} at 358 629
\pinlabel {${\cal K}^{(1)}$} <0pt,4pt> at 199 689
\pinlabel {${\cal K}^{(2)}$} <0pt,4pt> at 199 668
\pinlabel {flype} [b] at 216 585
\endlabellist
\includegraphics [height=40mm]{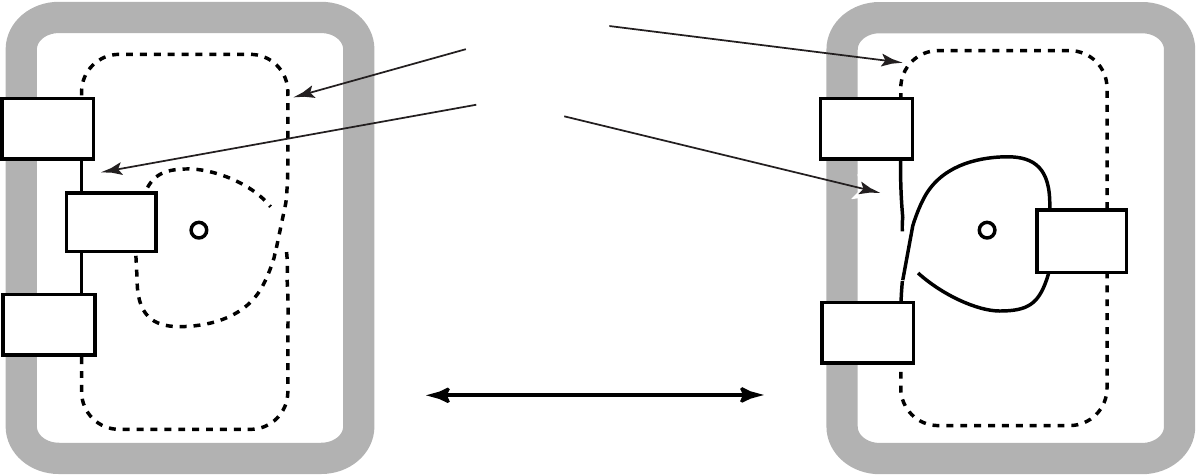}
\end{center}
\caption{A flype move}\label{flype}
\end{figure}
The thick gray arcs are parallel braid strands of ${\cal L}.$
Braidings occur inside the boxes $P, Q, R.$ In particular, box $R$ contains
even number of half twists of  ${\cal K}^{(1)}$ (dashed arc) and
${\cal K}^{(2)}$ (black
arc).
The flype move preserves the
number of braid strands and the algebraic crossing number of the link,
but it changes the algebraic crossing numbers of
link components ${\cal K}^{(1)}, {\cal K}^{(2)}.$
Namely, in the passage from the left sketch
to the right sketch, the algebraic crossing number of
${\cal K}^{(1)}$ decreases by
$1$ and the one for ${\cal K}^{(2)}$ increases by $1.$ (This means that
a flype move in general cannot be a composition of exchange moves.)}
\end{remark}

\begin{proof}[Proof of \fullref{cable-thm-link}] 
Suppose
$L=K^{(1)}\cup\cdots\cup K^{(l)}$ is a minimal braid
representative of ${\cal L}={\cal K}^{(1)}\cup\cdots\cup{\cal
K}^{(l)},$ ie, $b_L=b_{\cal L}.$ Let $c_{x,y}:=2\cdot{\rm
lk}({\cal K}^{(x)},{\cal K}^{(y)}) \mbox{ for } x\neq y.$ Then
\begin{equation}\label{c_L}
c_{\cal L}=c_L=\sum_{1\leq x <y\leq l}c_{x,y} + \sum_{j=1}^lc_{{\cal K}^{(j)}}.
\end{equation}
Let $k_j, n_j, 0\leq m_j<p$ be integers with
\begin{equation}\label{q_j}
q_j=p\cdot c_{{\cal K}^{(j)}}+k_j=p(c_{{\cal K}^{(j)}}+n_j)+m_j \quad \mbox{ for } j=1,\cdots,l.
\end{equation}

Williams proved that the braid index of ${\cal L}'={\cal
K}^{(1)}_{p,q_1}\cup \cdots\cup{\cal K}^{(l)}_{p,q_l}$ is $p\cdot
b_{\cal L}$ \cite{W}. Let $L'$ be a minimal braid representative
of ${\cal L}'.$ Let $N_j$ be a tubular neighborhood of $K^{(j)}.$
Let $\phi$ be a compact support diffeomorphism of $\mathbb{R}^3$ such that
$\phi(L')=:{K'}^{(1)}\cup\cdots\cup {K'}^{(l)} \subset
\partial\phi(N_1\cup\cdots N_l)$ has a minimal braid position.

\newpage
As we did in \eqref{K_i}, for each $j=1,\cdots,l,$
construct $p$ parallels $K^{(j)}_1,\cdots,K^{(j)}_p \subset \partial\phi(N_j)$
from each ${K'}^{(j)}$
by cutting out an inner most meridian disk $d_j\subset \phi(N_j)$
and adding an $m_j/p$--twist along annulus $A_j \subset \partial\phi(N_j)$
then gluing. Thus,
\begin{equation}\label{lk}
{\rm lk}(K_i^{(j)},K^{(h)})= \left\{
\begin{array}{ll}
c_{{\cal K}^{(j)}}+n_j & \ \mbox{when } j=h, \\
{\rm lk}(K^{(j)},K^{(h)})= \frac{1}{2}c_{j,h} & \ \mbox{otherwise. }
\end{array}
\right.
\end{equation}

Let $$L_i:=K^{(1)}_i \cup\cdots\cup =K^{(l)}_i \ \mbox{ for } i=1,\cdots,p.$$
Thanks to \cite{W} we know that $L_i \simeq L$ and $b_{L_i}=b_L=b_{\cal L}.$
By assumption of \fullref{cable-thm-link}, it follows that
$c_{L_i}=c_{L}=c_{\cal L}.$
The $(p\cdot l)$--component link $L_1 \cup\cdots\cup L_p$ has the algebraic crossing
number;
\begin{eqnarray*}
c_{L_1 \cup\cdots\cup L_p} &=& \sum_{i=1}^p c_{L_i} +
                               \sum_{x\neq y}{\rm lk}(L_x, L_y) \\
                           &=& p\cdot c_{\cal L} +\sum_{i=1}^p(p-1){\rm lk}(L_i,L)\\
                           &=& p\cdot c_{\cal L} +\sum_{i=1}^p(p-1)
                               \{ \sum_{j=1}^l{\rm lk}(K_i^{(j)},K^{(j)})
                                +\sum_{x\neq y}{\rm lk}(K_i^{(x)},K^{(y)})\} \\
                           &\stackrel{\eqref{lk}}{=}& p\cdot c_{\cal L}
          + p(p-1)\{\sum_{j=1}^l(c_{{\cal K}^{(j)}}+n_j)+\sum_{x<y}c_{x,y}\} \\
         & \stackrel{\eqref{c_L}}{=}& p\cdot c_{\cal L}
                   + p(p-1)(c_{\cal L}+\sum_{j=1}^ln_j) \\
         &=& p^2 c_{\cal L}+ p(p-1)(\sum_{j=1}^ln_j).
\end{eqnarray*}

Since only the difference between $L_1 \cup\cdots\cup L_p$ and $L'$ occurs on
the annuli $A_1,\cdots,A_l,$ we have
\begin{eqnarray*}
c_{L'}&=&c_{L_1 \cup\cdots\cup L_p}+(p-1)\sum_{j=1}^lm_j \\
      &=& p^2 c_{\cal L}+(p-1)\sum_{j=1}^l(pn_j+m_j)     \\
      &\stackrel{\eqref{q_j}}{=}& p^2 c_{\cal L}+(p-1)\sum_{j=1}^lk_j,
\end{eqnarray*}
which is independent of the choice of braid representative $L'$.
\end{proof}

With regard to \fullref{greedy} we have:

\begin{theorem}\label{cable-greedy-conj}
Let ${\cal L}={\cal K}^{(1)}\cup\cdots\cup{\cal K}^{(l)}$ be an $l$--component
link satisfying all the assumptions in \fullref{cable-thm-link}.
If \fullref{greedy} is true for ${\cal L}$ then
it is also true for its $p$--cable
${\cal L}':={\cal K}^{(1)}_{p,q_1}\cup\cdots\cup{\cal K}^{(l)}_{p,q_l}.$
\end{theorem}

\begin{proof}[Proof of \fullref{cable-greedy-conj}] 
Let $L$ (resp.\ $L'$) be a braid representative of ${\cal L}$ (resp.\
${\cal L}'$). Take tubular neighborhoods $N=N_1\cup\cdots\cup N_l$
of $L$ (each $N_j$ is a solid torus) and let
$\phi\co \mathbb{R}^3 \to\mathbb{R}^3$ be a compact support
diffeo morphism such that
$\phi(L)=:K^{(1)}\cup\cdots\cup K^{(l)}$,
$\phi(L')=:K'^{(1)}\cup\cdots\cup K'^{(l)}$ have braid positions.
They are not necessarily minimal braid representatives and in
general $b_{\phi(L')}\neq p\cdot b_{\phi(L).}$ We may assume that
$K'^{(j)}\subset\partial\phi(N_j)\simeq T^2.$

Let plane $\pi_j=\{(x,y,z)|y=y_0\},$ innermost meridian loop
$J_j\subset \pi_j\cap\partial\phi(N_j),$ and meridian disk
$d_j\subset \pi_j\cap\phi(N_j)$ be as in the proof of \fullref{cable-thm}. We may assume that the braid axis is {\em not}
contained in $\pi_j$ ie, $y_0\neq 0.$

We deform $K'^{(j)}$ in the following way: Suppose sub-arcs
$u\subset K'^{(j)}$ and $v\subset J_j$ bound a disk ${\cal
D}\subset\partial\phi(N_j)\simeq T^2$. If ${\cal D}$ is innermost,
then replace $u$ with $v$. Repeat this until $K'^{(j)}$ and $J_j$
do not bound any disk in $T^2.$ Add up all the linking numbers of
$\partial{\cal D}$'s with the $z$--axis and denote it by
$x^{(j)}\geq 0.$

Next, from the deformed $K'^{(j)}$ above, construct $p$--parallels
$K'^{(j)}_1,\cdots,K'^{(j)}_p$ as in \eqref{K_i}. Since the plane
$\pi_j$ does not contain the $z$--axis, the $m_j/p$--twist along a
thin annulus does not change the number of braid strands.

Suppose that $b_{K^{(j)}}=b_{{\cal K}^{(j)}}+y^{(j)}$ and
$b_{K'^{(j)}_i}=b_{{\cal K}^{(j)}}+y^{(j)}+z^{(j)}_i$ with
$y^{(j)},\ y^{(j)}+z^{(j)}_i\geq 0.$ Let $L_i:=K^{(1)}_i
\cup\cdots\cup =K^{(l)}_i$ then
\begin{equation}\label{b_L_i}
b_{L_i} =\sum_{j=1}^l b_{K^{(j)}_i}
        = b_{\cal L} + \sum_{j=1}^l (y^{(j)}+z^{(j)}_i),
\end{equation}
\begin{equation}\label{b_Phi(L')}
b_{\phi(L')}=\sum_{i=1}^p b_{L_i}+ \sum_{j=1}^lx^{(j)}
        \stackrel{\eqref{b_L_i}}{=} p\cdot b_{\cal L} +
        \sum_{j=1}^l(x^{(j)}+p\cdot y^{(j)}+\sum_{i=1}^pz^{(j)}_i).
\end{equation}
Since $L_i\simeq {\cal L},$ our assumption of this theorem and \eqref{b_L_i}
give us
\begin{equation}\label{c_L_i}
c_{\cal L}-\sum_{j=1}^l (y^{(j)}+z^{(j)}_i) \leq c_{L_i} \leq
c_{\cal L}+\sum_{j=1}^l (y^{(j)}+z^{(j)}_i).
\end{equation}
As in the proof of \fullref{cable-thm-link}, let $k_j, n_j,
0\leq m_j<p$ satisfy $ q_j=p\cdot c_{{\cal
K}^{(j)}}+k_j=p(c_{{\cal K}^{(j)}}+n_j)+m_j$. Then we have
\begin{eqnarray*}
c_{L_1 \cup\cdots\cup L_p} &=& \sum_{i=1}^p c_{L_i} +
                               \sum_{x\neq y}{\rm lk}(L_x, L_y) \\
                           &\stackrel{\eqref{c_L_i}}{\leq}&
                           p(c_{\cal L}+\sum_{j=1}^l y^{(j)})+
                           \sum_{i=1}^p \sum_{j=1}^lz^{(j)}_i +
                           p(p-1)(c_{\cal L}+\sum_{j=1}^ln_j) \\
                           &=& p^2 c_{\cal L}+ p(p-1)(\sum_{j=1}^ln_j) +
                           \sum_{j=1}^l (p\cdot y^{(j)}+
                           \sum_{i=1}^p z^{(j)}_i),
\end{eqnarray*}
and
\begin{eqnarray}
c_{\phi(L')} &\leq& c_{L_1 \cup\cdots\cup L_p} +\sum_{j=1}^l m_j +
                    \sum_{j=1}^lx^{(j)}    \nonumber               \\
             &\leq& p^2 c_{\cal L}+ (p-1)\sum_{j=1}^l k_j +
                    \sum_{j=1}^l(x^{(j)}+p\cdot y^{(j)}+
                    \sum_{i=1}^pz^{(j)}_i). \label{upper}
\end{eqnarray}
Similarly,
\begin{equation}\label{lower}
p^2 c_{\cal L}+ (p-1)\sum_{j=1}^l k_j -
\sum_{j=1}^l(x^{(j)}+p\cdot y^{(j)}+\sum_{i=1}^pz^{(j)}_i) \leq
c_{\phi(L')}.
\end{equation}
We conclude the theorem by \eqref{cable-thm-link-eq},
\eqref{b_Phi(L')}, \eqref{upper} and \eqref{lower}. 
\end{proof}

\begin{corollary}\label{iterated-torus}
Conjectures \ref{Jones-conj} and \ref{greedy} apply to iterated torus knots.
\end{corollary}

\begin{proof}[Proof of \fullref{iterated-torus}]
We know that the both conjectures apply to torus knots (\fullref{ex-of-conj}). Thanks to Theorems \ref{cable-thm}, \ref{cable-thm-link} 
and \ref{cable-greedy-conj}, we have this corollary.
\end{proof}

\subsection{Connect sum and the conjecture}
We will prove the following:

\begin{theorem}\label{sum}
If \fullref{Jones-conj} is true for knot types ${\cal K}^1$ and
${\cal K}^2$ then it is also true for the connect sum ${\cal K}^1\sharp{\cal K}^2$.

In particular, denoting the unique algebraic crossing numbers of ${\cal K}^i$
in minimal braid representatives by $c_{{\cal K}^i}$ we have
$$c_{{\cal K}^1\sharp {\cal K}^2} = c_{{\cal K}^1} + c_{{\cal K}^2}.$$
\end{theorem}

Before we prove \fullref{sum} let us recall two important known results:

\begin{lemma}\label{unique-prime}{\rm \cite[Theorem 2.12]{Lickorish}}\qua
Up to ordering of summands, there is a unique expression for a
knot type ${\cal K}$ as a finite connect sum of prime knots.
\end{lemma}

\begin{figure}[ht!]
\begin{center}
\labellist\small
\pinlabel {$P$} at 240 347
\pinlabel {$Q$} at 241 241
\endlabellist
\includegraphics [height=30mm]{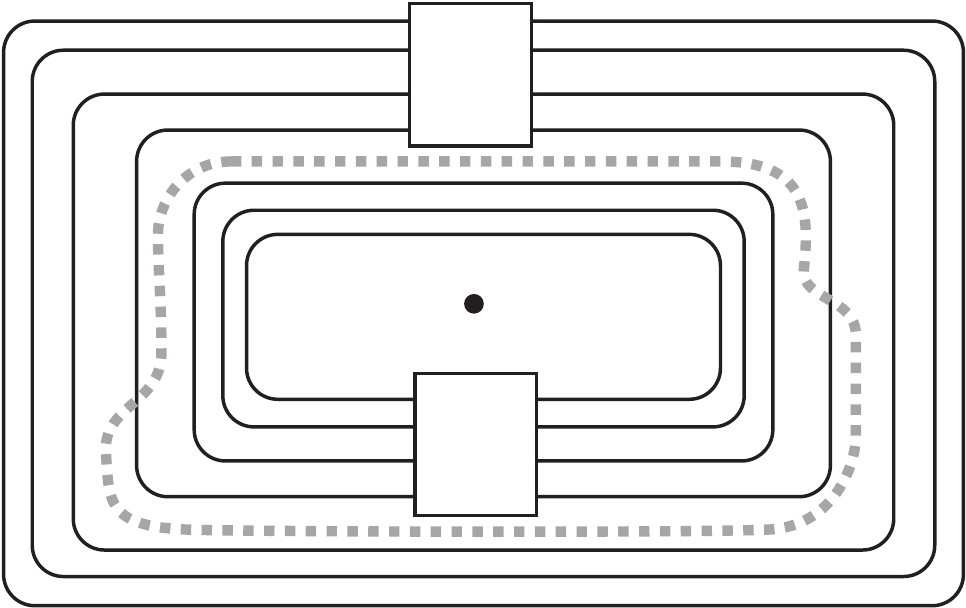}
\end{center}
\caption{An obvious composite braid} \label{obvious-sum}
\end{figure}

\begin{lemma}\label{comp-th}
{\rm (\textbf{The composite braid theorem,} \cite{BM4}.)} Let
${\cal K}$ be a composite link, and let $K$ be an arbitrary closed
$n$--braid representative of ${\cal K}$.  Then there is an obvious composite
$n$--braid representative $K^\bullet$ of ${\cal K}$ (see \fullref{obvious-sum}) and a finite
sequence of closed $n$--braids:
$$K=K_0 \rightarrow K_1 \rightarrow \cdots \rightarrow K_m=K^\bullet$$
such that $K_{i+1}$ is obtained from $K_i$ by either braid isotopy
or an exchange move.
\end{lemma}

\begin{proof}[Proof of \fullref{sum}] 
Since an exchange move
does not change the algebraic crossing number, Lemmas \ref{unique-prime} and
\ref{comp-th} imply the truth of \fullref{sum}.
\end{proof}

As a corollary of \fullref{sum} we have:

\begin{theorem}\label{sum'}
If \fullref{greedy} is true for ${\cal K}^1, {\cal K}^2$ then it
is also true for ${\cal K}^1 \sharp {\cal K}^2.$
\end{theorem}

\begin{proof}[Proof of \fullref{sum'}] 
Let $K$ be a braid
representative of ${\cal K}^1 \sharp {\cal K}^2.$ By \fullref{comp-th}, after applying exchange moves and braid isotopy to
$K$ one can get a composite braid representative $K^\bullet=K^1
\sharp K^2.$ Suppose $b_{K^i}=b_{{\cal K}^i}+x_i$ with $x_i\geq
0.$ Then
\begin{eqnarray}
b_K &=& b_{K^\bullet}=b_{K^1}+b_{K^2}-1 =(b_{{\cal K}^1}+b_{{\cal
          K}^2}-1)+(x_1+x_2) \nonumber \\
    &=& b_{{\cal K}^1 \sharp {\cal K}^2} +(x_1+x_2). \label{b_K}
\end{eqnarray}
Our assumption gives $c_{{\cal K}^i}-x_i\leq c_{K^i}\leq c_{{\cal
K}^i}+x_i.$ Since $c_K=c_{K^\bullet}=c_{K^1}+c_{K^2},$ we have
$(c_{{\cal K}^1}+c_{{\cal K}^2})-(x_1+x_2)\leq c_K \leq (c_{{\cal
K}^1}+c_{{\cal K}^2})+(x_1+x_2).$ Thanks to \fullref{sum},
\begin{equation}\label{c_K}
c_{{\cal K}^1 \sharp {\cal K}^2}-(x_1+x_2)\leq c_K
\leq c_{{\cal K}^1 \sharp {\cal K}^2}+(x_1+x_2).
\end{equation}
The truth of \fullref{greedy} follows by \eqref{b_K}, \eqref{c_K}.
\end{proof}

\bibliographystyle{gtart}
\bibliography{link}

\begin{thebibliography}{}
\providecommand\bibmarginpar{\leavevmode\marginpar}
\def\urlstyle#1{{\tt #1}}

\bibitem{Ben}
\textbf{D Bennequin}, \emph{Entrelacements et \'equations de {P}faff}, from:
  ``Third Schnepfenried geometry conference, Vol. 1 (Schnepfenried, 1982)'',
  Ast\'erisque 107, Soc. Math. France, Paris (1983)  87--161 \xox{MR}{753131}

\bibitem{B}
\textbf{J\,S Birman}, \emph{Braids, links, and mapping class groups}, Annals of
  Mathematics Studies 82, Princeton University Press (1974) \xox{MR}{0375281}

\bibitem{BM4}
\textbf{J\,S Birman}, \textbf{W\,W Menasco},
  \href{http://dx.doi.org/10.1007/BF01233423} {\emph{Studying links via closed
  braids. {IV}. {C}omposite links and split links}}, Invent. Math. 102 (1990)
  115--139 \xox{MR}{1069243}

\bibitem{BM3}
\textbf{J\,S Birman}, \textbf{W\,W Menasco},
  \href{http://projecteuclid.org/getRecord?id=euclid.pjm/1102623463}
  {\emph{Studying links via closed braids. {III}. {C}lassifying links which are
  closed {$3$}-braids}}, Pacific J. Math. 161 (1993) 25--113 \xox{MR}{1237139}

\bibitem{MTWS-I}
\textbf{J\,S Birman}, \textbf{W\,W Menasco},
  \href{http://dx.doi.org/10.2140/gt.2006.10.413} {\emph{Stabilization in the
  braid groups. {I}. {MTWS}}}, Geom. Topol. 10 (2006) 413--540
  \xox{MR}{2224463}

\bibitem{FW}
\textbf{J Franks}, \textbf{R\,F Williams},
  \href{http://links.jstor.org/sici?sici=0002-9947(198709)303:1%3C97:BATJP%3E2%
.0.CO%3B2-N} {\emph{Braids and the {J}ones polynomial}}, Trans. Amer. Math.
  Soc. 303 (1987) 97--108 \xox{MR}{896009}

\bibitem{Giroux}
\textbf{E Giroux}, \emph{G\'eom\'etrie de contact: de la dimension trois vers
  les dimensions sup\'erieures}, from: ``Proceedings of the International
  Congress of Mathematicians, Vol. II (Beijing, 2002)'', Higher Ed. Press,
  Beijing (2002)  405--414 \xox{MR}{1957051}

\bibitem{Jones-1}
\textbf{V\,F\,R Jones},
  \href{http://links.jstor.org/sici?sici=0003-486X(198709)2:126:2%3C335:HAROB%
G%3E2.0.CO%3B2-E} {\emph{Hecke algebra representations of braid groups and link
  polynomials}}, Ann. of Math. $(2)$ 126 (1987) 335--388 \xox{MR}{908150}

\bibitem{K}
\textbf{K Kawamuro}, \emph{Conjectures on the braid index and the algebraic
  crossing number}, to appear in proceedings of ``Intelligence of Low
  Dimensional Topology (Hiroshima, 2006)''

\bibitem{Lickorish}
\textbf{W\,B\,R Lickorish}, \emph{An introduction to knot theory}, Graduate
  Texts in Mathematics 175, Springer (1997) \xox{MR}{1472978}

\bibitem{MM}
\textbf{P\,M Melvin}, \textbf{H\,R Morton}, \emph{Fibred knots of genus {$2$}
  formed by plumbing {H}opf bands}, J. London Math. Soc. $(2)$ 34 (1986)
  159--168 \xox{MR}{859157}

\bibitem{Morton}
\textbf{H\,R Morton}, \emph{Seifert circles and knot polynomials}, Math. Proc.
  Cambridge Philos. Soc. 99 (1986) 107--109 \xox{MR}{809504}

\bibitem{MS}
\textbf{H\,R Morton}, \textbf{H\,B Short}, \emph{The {$2$}-variable polynomial
  of cable knots}, Math. Proc. Cambridge Philos. Soc. 101 (1987) 267--278
  \xox{MR}{870598}

\bibitem{Murasugi}
\textbf{K Murasugi},
  \href{http://links.jstor.org/sici?sici=0002-9947(199107)326:1%3C237:OTBIOA%3%
E2.0.CO%3B2-R} {\emph{On the braid index of alternating links}}, Trans. Amer.
  Math. Soc. 326 (1991) 237--260 \xox{MR}{1000333}

\bibitem{N}
\textbf{W\,D Neumann}, \emph{Private communication}

\bibitem{NR}
\textbf{W\,D Neumann}, \textbf{L Rudolph},
  \href{http://dx.doi.org/10.1016/0040-9383(90)90026-G} {\emph{Difference index
  of vectorfields and the enhanced {M}ilnor number}}, Topology 29 (1990)
  83--100 \xox{MR}{1046626}

\bibitem{NZ}
\textbf{W\,D Neumann}, \textbf{D Zagier},
  \href{http://dx.doi.org/10.1016/0040-9383(85)90004-7} {\emph{Volumes of
  hyperbolic three-manifolds}}, Topology 24 (1985) 307--332 \xox{MR}{815482}

\bibitem{S-preprint}
\textbf{A Stoimenow}, \emph{Properties of closed $3$-braids}
  \xox{arXiv}{math.GT/0606435}

\bibitem{S}
\textbf{A Stoimenow}, \href{http://dx.doi.org/10.1090/S0002-9947-02-03022-2}
  {\emph{On the crossing number of positive knots and braids and braid index
  criteria of {J}ones and {M}orton-{W}illiams-{F}ranks}}, Trans. Amer. Math.
  Soc. 354 (2002) 3927--3954 \xox{MR}{1926860}

\bibitem{W}
\textbf{R\,F Williams},
  \href{http://projecteuclid.org/getRecord?id=euclid.pjm/1102635274} {\emph{The
  braid index of generalized cables}}, Pacific J. Math. 155 (1992) 369--375
  \xox{MR}{1178031}

\bibitem{Xu}
\textbf{P Xu}, \href{http://dx.doi.org/10.1142/S0218216592000185} {\emph{The
  genus of closed {$3$}-braids}}, J. Knot Theory Ramifications 1 (1992)
  303--326 \xox{MR}{1180404}

\end{thebibliography}

\end{document}